\documentclass[hidelinks]{siamart251216}

\usepackage[utf8]{inputenc}
\usepackage[T1]{fontenc}
\usepackage{mathtools,amssymb,mathrsfs,bm}
\usepackage{microtype}
\newcommand{\IR}{{\mathbb R}}
\newcommand{\IC}{{\mathbb C}}

\newcommand{\frakT}{{\mathfrak T}}
\newcommand{\D}{\text{\normalfont D}}
\newcommand{\x}{{\bf{x}}}
\newcommand{\y}{{\bf{y}}}
\newcommand{\br}{{\normalfont{\boldsymbol{r}}}}
\newcommand{\bj}{\bm{j}}
\newcommand{\bv}{\bm{v}}
\newcommand{\bu}{\bm{u}}
\newcommand{\bw}{\bm{w}}

\newcommand{\bn}{\bm{n}}
\newcommand{\bp}{{\bm p}}
\newcommand{\bE}{\bm{E}}
\newcommand{\bV}{\bm{V}}
\newcommand{\bL}{\boldsymbol{L}}
\newcommand{\bH}{\boldsymbol{H}}
\newcommand{\bff}{\bm{f}}
\newcommand{\bphi}{\boldsymbol \phi}
\newcommand{\bpsi}{\boldsymbol \psi}
\newcommand{\bseta}{{\boldsymbol \eta}}

\newcommand{\btau}{{\boldsymbol \tau}}
\newcommand{\einc}{\bE^{\mathrm{inc}}}
\newcommand{\escat}{\bE^{\mathrm{scat}}}
\newcommand{\half}{\frac{1}{2}}

\newcommand{\widehatGamma}{\widehat\Gamma}
\newcommand{\widehatD}{\widehat{\D}}
\newcommand{\widehatX}{\widehat X}
\newcommand{\widehatY}{\widehat Y}

\newcommand{\cA}{\mathcal A}

\newcommand{\cE}{\mathcal E}

\newcommand{\cO}{\mathcal O}
\newcommand{\R}{{\mathcal R}}

\newcommand{\V}{{\mathsf V}}
\newcommand{\dd}{\,{\rm d}}
\newcommand{\ds}{\dd s}
\newcommand{\norm}[1]{\left \lVert #1 \right \rVert}
\newcommand{\snorm}[1]{\left \lvert #1 \right \rvert}
\newcommand{\dual}[2]{\left\langle#1,#2\right\rangle}
\renewcommand{\div}{\operatorname{div}}
\newcommand{\divGh}{\div_{\widehat\Gamma}}
\newcommand{\curl}{\operatorname{curl}}
\newcommand{\curlPiola}[1]{\mathsf T_{#1}^{\operatorname{curl}}}

\newsiamthm{assumption}{Assumption}
\newsiamremark{remark}{Remark}

\headers{Shape Holomorphy of the EFIE on Lipschitz Surfaces}{F. Henr\'iquez}
\title{Shape Holomorphy of the Maxwell Electric Field Integral Equation on Lipschitz Surfaces\thanks{\textbf{Funding:} FH's work is funded by the Deutsche Forschungsgemeinschaft (DFG, German
Research Foundation), Project-ID 258734477, SFB 1173, and the Austrian
Science Fund (FWF) under the project I6667-N. Funding was also received from
the European Research Council (ERC) under the Horizon
Europe research and innovation programme of the European Union (Grant agreement No.~101125225,
Project Name: ``New Frontiers in Optimal Adaptivity'', PI: Prof.\ Dr.\
Michael Feischl).}}
\author{Fernando Henr\'iquez\thanks{Institute for Analysis and Scientific Computing,
Technische Universit\"{a}t Wien, Wiedner Hauptstrasse 8-10, A-1040 Wien, Austria
(\email{fernando.henriquez@asc.tuwien.ac.at}).}}
\hypersetup{
 pdftitle={Shape Holomorphy of the Maxwell Electric Field Integral Equation on Lipschitz Surfaces},
 pdfauthor={Fernando Henriquez}
}

\begin{document}
\maketitle

\begin{abstract}
We establish the holomorphic dependence of the Maxwell electric field
integral equation (EFIE) operator upon the shape of a family of bounded, perfectly
conducting obstacles with connected Lipschitz boundaries, a property
referred to as \emph{shape holomorphy}. More precisely, after
transporting tangential fields and the EFIE operator to a fixed reference
boundary, we prove that the map from the shape deformation to the
transported operator, the \emph{shape-to-operator map}, admits a bounded holomorphic extension to a complex
neighborhood of a compact set of admissible bijective and bi-Lipschitz shape transformations.

The EFIE operator acts on the space $\bH^{-\half}(\div_{\Gamma},\Gamma)$
of tangential fields, which itself depends on the boundary $\Gamma$.
The main technical step is therefore the transport of tangential fields
between Lipschitz boundaries while preserving this property.
In this work, this is carried out by lifting tangential traces to an
$\bH(\curl)$ field, transporting them with a curl-conforming volume
pullback, and, finally, taking tangential traces on the reference boundary. These
three steps give the contravariant surface Piola pullback for the surface
current, and the covariant pullback is obtained analogously using the
tangential projection trace and its lifting.
Using these operators, the transported EFIE operator is
expressed in terms of the scalar Helmholtz single layer operator
transported to the reference boundary, for which the shape holomorphy
on Lipschitz surfaces has already been established, and standard surface differential operators.

The presently established result extends previous work concerning the EFIE on closed
$\mathscr C^{1,1}$ surfaces and on the Laplace and Helmholtz Calder\'on
projectors on Lipschitz surfaces.
Furthermore, it provides the analytic foundation for uncertainty quantification and model order reduction in time-harmonic electromagnetic scattering on families
of parametrically varying Lipschitz boundaries, and the theoretical tools
to establish shape holomorphy of the complete Maxwell Calder\'on
projector.
\end{abstract}



\raggedbottom
\section{Introduction}

In this work, we consider the electric field integral equation (EFIE)
for time-harmonic electromagnetic scattering by bounded, perfectly
conducting obstacles in $\IR^3$, each with a connected
Lipschitz boundary.
The scattered electric field satisfies Maxwell's equations in the
exterior of the obstacle and the Silver--M\"uller radiation condition
at infinity, while the tangential component of the total electric
field vanishes on the boundary. We study the holomorphic dependence
of the EFIE operator and its solution on shape deformations after
transport to a fixed reference boundary,
a property referred to as \emph{shape holomorphy}.

Constructing efficient surrogates for the EFIE operator and the
corresponding boundary density is a central task in forward and inverse
uncertainty quantification (UQ) and reduced order modelling of
electromagnetic scattering problems. In these applications, the
surrogates must capture how the operator and its solution vary across
a family of obstacle geometries. Describing this family through
geometric parameters gives rise to the parameter-to-operator and
parameter-to-solution maps to be approximated. Such descriptions may
involve many, or even countably infinitely many, parameters, making
the approximation of these maps computationally challenging. For
tensor product methods, even a fixed resolution in each variable leads
to a cost that grows exponentially with the number of retained
parameters, a manifestation of the \emph{curse of dimensionality}.
This motivates approximation methods that exploit regularity with
respect to the geometric parameters. Establishing this regularity
for the EFIE operator and its solution is therefore a key step in
the design and analysis of these methods.

In \cite{CCS15}, it has been established that the holomorphic
dependence on the parametric input, together with anisotropic
parameter dependence and suitable summability assumptions, enables
sparse polynomial approximation with algebraic convergence rates
independent of the number of parameters. This regularity is quantified
through bounded holomorphic extensions to suitable complex domains,
whose sizes reflect the dependence on the individual parameters.
A variety of computational methods draw on this regularity in their
convergence analysis, including sparse grid methods
\cite{ZS17,SS13,NTW2008,HHPS2018} and higher-order quasi-Monte Carlo
integration \cite{DKL14,DLC16}. The same regularity also underlies
approximation estimates for neural network surrogates
\cite{SZ19,DH2026,HS26NeuralShape}.

\subsection{Previous Work}
\label{sec:previous-work}

The holomorphic dependence of solution maps on shape deformations has
been studied for a variety of boundary value problems. For variational
formulations posed in the domain, holomorphy with respect to the
parameters describing the geometry has been established for elliptic
diffusion \cite{CNT2016,HPS2016,HS2022} and acoustic scattering
\cite{HSSS2018}. Holomorphy with respect to the domain transformation
itself is established for stationary Navier--Stokes flow
\cite{CSZ18} and electromagnetic scattering by Lipschitz obstacles
\cite{JSZ2017}. In the latter work, a curl-conforming pullback
transports the electric field to a fixed reference $\bH(\curl)$
space, where its holomorphic dependence on the deformation is
established. For volume integral formulations, shape holomorphy
of the acoustic Lippmann--Schwinger operator itself is established
in \cite{henriquez2024domain} under Lipschitz deformations of the
inhomogeneous scatterer.

For boundary integral operators (BIOs), shape derivatives under real
boundary perturbations have been characterized in acoustic, elastic,
and electromagnetic scattering
\cite{Pot94,Charam95,Pot96,Pot96_3,LC121,CostabelLeLouer2012ShapeII}.
Such differentiability, even to all orders, does not by itself imply
the existence of a holomorphic extension. For scalar BIOs, shape
holomorphy has been established for the Laplace and Helmholtz
Calder\'on projectors on closed $\mathscr C^2$ curves
\cite{HS21,henriquez2021shape} and for weakly singular and hypersingular
operators on multiple open arcs \cite{PHJ23}. At Lipschitz regularity,
\cite{dolz2024parametric} establishes parametric holomorphy for weakly
singular operators under affine boundary deformations, with the
analysis carried out in an $L^2$ operator setting.

Recently, in \cite{HDPS26Laplace,HDPS26Helmholtz}, we established shape
holomorphy of the full Laplace and Helmholtz Calder\'on projectors on
Lipschitz boundaries. Our analysis expresses the four BIOs as
compositions of volume potentials, trace operators, and their
adjoints, and establishes holomorphy of these constituent operators
separately. Their composition then yields bounded holomorphic
extensions of the BIOs in their natural fractional Sobolev trace
spaces. In the Helmholtz case, we also obtain a common domain of
holomorphic extension for all positive wavenumbers, together with
operator bounds explicit in both the wavenumber and the size of
this domain of holomorphy.

Composition with $\br$ transports a continuous scalar field from the
deformed boundary $\Gamma_\br=\br(\widehatGamma)$ to the reference
boundary $\widehatGamma$. For a continuous tangential field, however,
the same composition may introduce a normal component on
$\widehatGamma$, since the tangent planes may differ at corresponding
points. Under small smooth deformations, orthogonal projection removes
this component and identifies the spaces of continuous tangential
fields on the two boundaries, as proposed in \cite{Pot96_3}. The EFIE trace space
imposes the additional requirement that the weak surface divergence
belong to $H^{-\half}(\Gamma_\br)$. As observed in
\cite{CostabelLeLouer2012ShapeII}, checking this condition for a
projected field requires transporting it back to $\Gamma_\br$.
Consequently, the image of
$\bH^{-\half}(\div_{\Gamma_\br},\Gamma_\br)$ under this projected
pullback still depends on the deformation.

For smooth, simply connected boundaries, the pullback proposed in
\cite{CostabelLeLouer2012ShapeII} uses the surface Helmholtz
decomposition given in \cite{DeLaBourdonnaye1993}, which expresses the tangential trace space
$\bH^{-\half}(\div_{\Gamma_\br},\Gamma_\br)$ as the direct sum
of surface gradients of $H^{\frac{3}{2}}(\Gamma_\br)$ potentials and
surface curls of $H^{\half}(\Gamma_\br)$ potentials. The two scalar potentials are pulled
back to $\widehatGamma$, where the surface gradient of the first and
the surface curl of the second reconstruct the tangential field.
Each Maxwell BIO maps the input potentials to the output potentials
through scalar BIOs and surface operators. Using \cite{LC121} for the shape derivatives of the scalar BIOs,
the authors combine these with those of the
surface differential operators and the inverse Laplace--Beltrami operator. This yields
G\^ateaux differentiability to all orders of the electric and magnetic
BIOs on the fixed Maxwell trace space without loss of regularity.

A different approach, based on the contravariant surface Piola
transformation, is developed in \cite{YL16} for closed
$\mathscr C^{1,\alpha}$ reference surfaces, $0<\alpha\leq 1$.
This pullback transports tangential fields while transforming their
weak surface divergence as a scalar density. The authors derive
explicit formulas for the first Fr\'echet derivatives of the
transported electric and magnetic BIOs with respect to the boundary
parametrization and state differentiability to all orders.

Using this Piola pullback, \cite{EscapilJerez2026EFIE} establishes
parametric holomorphy of the EFIE operator under countably
parametrized affine deformations of a closed $\mathscr C^{1,1}$
surface. Together with a summability assumption on the
$\mathscr C^{1,1}$ deformation fields, this regularity yields
dimension-independent convergence rates for sparse Legendre
approximation of the operator. Closed surfaces and the screens
they support are treated in \cite{Henriquez2026C11}, where we use
the scalar Helmholtz single layer to extend the EFIE operator
holomorphically to a complex neighborhood of the admissible boundary
transformations. 
We hasten to add that both holomorphy results are established under
$\mathscr C^{1,1}$ regularity assumptions on the closed surfaces
and their deformations.

\subsection{Contributions}
\label{sec:contributions}

We establish the holomorphic dependence of the EFIE operator on
admissible bi-Lipschitz deformations of perfectly conducting obstacles
in $\IR^3$ with connected Lipschitz boundaries. More precisely, after
transport to a fixed reference boundary, the operator admits a bounded
holomorphic extension in its natural Maxwell trace spaces
(Theorem~\ref{thm:shape-holomorphy}). For closed surfaces, this extends
the analyses of \cite{Henriquez2026C11,EscapilJerez2026EFIE} beyond
$\mathscr C^{1,1}$ regularity and permits Lipschitz domain
transformations. In particular, polyhedral reference domains are
included.

The passage to Lipschitz boundaries requires a transport that accounts
for both the tangential character of the boundary density and the weak
surface divergence appearing in the EFIE. We construct the required
boundary pullbacks through volume liftings. We lift a tangential field
on the deformed boundary to an $\bH(\curl)$ field in the obstacle,
apply the curl-conforming volume pullback, and take the corresponding
trace on the reference boundary. The resulting pullbacks preserve the
Maxwell duality pairing and provide a weak transformation rule for the
surface divergence. These properties allow the EFIE to be studied on
fixed Maxwell trace spaces at Lipschitz regularity. When the boundaries,
the deformation, and the tangential fields are smooth, the construction
recovers the classical contravariant surface Piola transformation for
the boundary density and the classical covariant surface Piola
transformation for the tangential electric field
(Remark~\ref{rem:classical-surface-piola}).

To establish holomorphy, we express the transported EFIE bilinear form
through the scalar Helmholtz single layer and two operator families on
the reference boundary. The transported intrinsic injection represents
the tangential boundary density in
$\bH^{-1/2}(\widehatGamma)$, allowing the single layer to act
componentwise, while the transported surface divergence supplies the
scalar density in the second term. We prove that both maps admit
bounded holomorphic extensions between fixed spaces
(Proposition~\ref{prop:coefficient-identification}). Combining these
extensions with the transported Sobolev duality pairing from
\cite{HDPS26Laplace} and the Helmholtz single layer analysis of
\cite{HDPS26Helmholtz} yields shape holomorphy of the EFIE operator.

We also quantify the dependence of the operator bounds on the
wavenumber and the size of the holomorphic extension. The extension
domain can be chosen independently of $\kappa>0$. For
$\kappa\geq\kappa_0>0$ and sufficiently small extension size $\delta>0$,
the norm of the transported EFIE operator is bounded by
$C(1+\kappa)^2$ on the real family and by
$C(1+\kappa)^2e^{C\kappa\delta}$ throughout the complex extension
domain, with $C$ independent of $\kappa$ and $\delta$. Taking
$\delta$ proportional to $\kappa^{-1}$ within the admissible range
keeps the exponential factor uniformly bounded and retains an
operator norm bound of order $(1+\kappa)^2$ throughout the
corresponding complex neighborhood.

At a fixed wavenumber for which the EFIE is well posed throughout the
admissible family, the inverse operator admits a bounded holomorphic
extension to a possibly smaller complex neighborhood
(Corollary~\ref{cor:efie-solution}). Applying this inverse to
transported data that admit a bounded holomorphic extension gives
the same regularity for the transported boundary density. We verify
the data assumption for plane-wave incidence, thereby obtaining shape
holomorphy of the corresponding solution map.

For the admissible affine families considered here, we further
establish $(\boldsymbol b,p,\varepsilon)$-holomorphy of the
parameter-to-operator map
(Corollary~\ref{cor:efie-parametric}). The weights $\boldsymbol b$
are the $W^{1,\infty}$ norms of the deformation fields and are assumed
to belong to $\ell^p(\mathbb N)$ for some $p\in(0,1)$. Under the same
well-posedness and data assumptions, the parameter-to-solution map
has this regularity as well, possibly with a smaller $\varepsilon$.
We verify the geometric assumptions for sufficiently small affine
perturbations of a polyhedral reference domain. These results supply
the analytic regularity needed to apply sparse polynomial
approximation theory \cite{CCS15} to EFIE operators and boundary
densities over families of Lipschitz geometries.

\subsection{Outline}
\label{sec:outline}

Section~\ref{sec:preliminaries} recalls the Maxwell trace spaces
and constructs the volume pullbacks and surface Piola transformations.
Section~\ref{sec:transport} represents the density and its surface
divergence by holomorphic operator families on the reference boundary.
Section~\ref{sec:main-results} establishes shape holomorphy of the EFIE
operator and the corresponding solution map.
Section~\ref{sec:efie-parametric} establishes
parametric holomorphy for affine domain transformations.
Section~\ref{sec:concluding-remarks} provides concluding remarks.

\subsection{Notation}
\label{sec:notation}

For (complex) Banach spaces $E$ and $F$, we denote by
$\mathscr{L}(E,F)$ the space of bounded linear operators from $E$ into
$F$ and by $\mathscr{L}_{\mathrm{iso}}(E,F)$ the subset of isomorphisms,
i.e.,~bounded linear operators with a bounded inverse.
Recall that $\mathscr{L}(E,F)$ is a (complex) Banach space equipped
with the standard operator norm.
We also denote by $E'$ the dual space of $E$, i.e.,
$E'\coloneqq\mathscr{L}(E,\IC)$.
Dual pairings are denoted by $\dual{\cdot}{\cdot}_{E'\times E}$,
with $\dual{f}{u}_{E'\times E}=f(u)$.
This pairing and the Sobolev and Maxwell duality pairings are complex bilinear.
We denote by $\mathscr{L}_{\mathrm{bl}}(E\times F,\IC)$ the space of
bounded bilinear maps from $E\times F$ into $\IC$, equipped with the
standard norm.

For a given $T\in\mathscr{L}(E,F)$, we define its dual
$T'\in\mathscr{L}(F',E')$ as the unique operator satisfying
\begin{equation*}
	\dual{f}{Tu}_{F'\times F}=\dual{T'f}{u}_{E'\times E},
	\quad u\in E,\ f\in F'.
\end{equation*}

All function and distribution spaces are understood over $\IC$, unless
stated otherwise. For a domain or boundary $G$, denote by $L^2(G)$ and
$L^\infty(G)$ the usual Lebesgue spaces on $G$, and by $H^s(G)$ the
standard Sobolev spaces on $G$.
We use boldface for vector-valued spaces, with
$\bH^s(G)\coloneqq H^s(G;\IC^3)$.
On a Lipschitz boundary $G$, we identify $H^{-\half}(G)$ with the dual
space of $H^{\half}(G)$ in the $L^2(G)$ duality pairing.
The complex Hilbert spaces used below are complexifications of real
Hilbert spaces. On $H^{\half}(G)$, we use the complex-linear extension of the Riesz map
for the real Slobodeckij inner product. This is an isometric isomorphism
onto $H^{-\half}(G)$, obtained by conjugating the argument before
applying the usual Riesz map to the linear dual.
Hats indicate quantities on the reference domain or its boundary. We write $\widehat{\nabla}$ for the gradient on $\widehatD$.

The dot denotes the complex-bilinear Euclidean product.
For $\x\in\IR^3$ and $R>0$, we write $B(\x,R)$ for the open Euclidean
ball centered at $\x$. We use the Euclidean vector norm and its induced
matrix norm.
We equip $W^{1,\infty}(\D;\IC^3)$ with the norm
\begin{equation*}
	\norm{\bu}_{W^{1,\infty}(\D;\IC^3)}
	=\norm{\bu}_{L^\infty(\D;\IC^3)}
	+\norm{D\bu}_{L^\infty(\D;\IC^{3\times 3})},
\end{equation*}
where $D\bu$ denotes the Jacobian matrix. The corresponding real-valued
space is equipped with the same norm.

\flushbottom
\section{Preliminaries}
\label{sec:preliminaries}
We recall important definitions and concepts used extensively in this work.
Firstly, Subsection~\ref{sec:maxwell-traces} recalls the Maxwell trace spaces
and their duality. In Subsection~\ref{sec:geometry}, we introduce the
admissible domain transformations. In Subsection~\ref{sec:pullback-operators}, we
define the volume and boundary pullbacks. Finally, Subsection~\ref{sec:piola}
constructs the surface Piola transformations and establishes their main properties.

\subsection{Maxwell Domain and Trace Spaces}
\label{sec:maxwell-traces}

Throughout this subsection we follow \cite{BCS2002}, to which we refer
for proofs and further details, see also
\cite{BuffaCostabelSchwab2002,BuffaHiptmair2003}.
Let $\D\subset\IR^3$ be a bounded Lipschitz domain with boundary $\Gamma \coloneqq \partial \D$
and outward unit normal $\bn$.  We define
\begin{equation*}
	\bH(\curl,\D)
	\coloneqq
	\left\{
	\bu\in L^2(\D;\IC^3):\curl\bu\in L^2(\D;\IC^3)
	\right\},
\end{equation*}
equipped with the norm
\begin{equation*}
	\norm{\bu}_{\bH(\curl,\D)}^2
	\coloneqq
	\norm{\bu}_{L^2(\D;\IC^3)}^2
	+
	\norm{\curl\bu}_{L^2(\D;\IC^3)}^2,
\end{equation*}
and we denote by $\bH_0(\curl,\D)$ the closure of
$\mathscr C_0^\infty(\D;\IC^3)$ in the $\bH(\curl,\D)$ norm.
Similarly, we define
\begin{equation*}
	\bH(\div,\D)
	\coloneqq
	\left\{
	\bu\in L^2(\D;\IC^3):\div\bu\in L^2(\D)
	\right\},
\end{equation*}
equipped with the norm
\begin{equation*}
	\norm{\bu}_{\bH(\div,\D)}^2
	\coloneqq
	\norm{\bu}_{L^2(\D;\IC^3)}^2
	+
	\norm{\div\bu}_{L^2(\D)}^2,
\end{equation*}
and we denote by $\bH_0(\div,\D)$ the closure of
$\mathscr C_0^\infty(\D;\IC^3)$ in the $\bH(\div,\D)$ norm.

Furthermore, we denote by $\gamma_\Gamma:H^1(\D)\to H^{\half}(\Gamma)$
the Dirichlet trace on $\Gamma$, which is a bounded and
surjective operator and admits a bounded linear right inverse denoted in the following
by $\mathsf Z_{\D}:H^{\half}(\Gamma)\to H^1(\D)$, see, e.g., \cite[Theorem~3.37]{mclean2000strongly}.
For the sake of simplicity, we use the same notation for the componentwise
action of $\gamma_\Gamma$ and $\mathsf Z_\D$ on vector-valued functions.

\subsubsection{Traces of \texorpdfstring{$\bH^1(\D)$}{H1(D)} Fields}
For $\bu\in\bH^1(\D)$, we define the tangential trace
$\gamma_{\times,\Gamma}$ and the tangential components trace
$\gamma_{\top,\Gamma}$ as
\begin{equation*}
	\gamma_{\times,\Gamma}\bu\coloneqq\bn\times\gamma_\Gamma\bu
	\quad\text{and}\quad
	\gamma_{\top,\Gamma}\bu\coloneqq\bn\times(\gamma_\Gamma\bu\times\bn).
\end{equation*}
Their ranges
\begin{equation*}
	V_\top\coloneqq\gamma_{\top,\Gamma}\bH^1(\D)
	\quad\text{and}\quad
	V_\times\coloneqq\gamma_{\times,\Gamma}\bH^1(\D)
\end{equation*}
define complex Hilbert spaces equipped with the quotient norms
\begin{equation*}
	\norm{\bseta}_{V_\top}
	\coloneqq\inf_{\substack{\bu\in\bH^1(\D)\\ \gamma_{\top,\Gamma}\bu=\bseta}}
	\norm{\bu}_{\bH^1(\D)}
	\quad\text{and}\quad
	\norm{\bj}_{V_\times}
	\coloneqq\inf_{\substack{\bu\in\bH^1(\D)\\ \gamma_{\times,\Gamma}\bu=\bj}}
	\norm{\bu}_{\bH^1(\D)}.
\end{equation*}

We denote by $V_\top'$ and $V_\times'$ their dual spaces with respect
to the pivot space of square integrable tangential fields
\begin{equation*}
	\bL_t^2(\Gamma)
	\coloneqq
	\left\{
	\bv\in L^2(\Gamma;\IC^3):\bv\cdot\bn=0\text{ a.e. on }\Gamma
	\right\}.
\end{equation*}

\subsubsection{Tangential Differential Operators}
In the following, we denote by $\nabla_\Gamma$ the surface gradient and by
$\mathbf{curl}_\Gamma\coloneqq-\bn\times\nabla_\Gamma$ the vector surface curl
which are defined using a chart localization construction,
see, e.g., \cite[Definition~3.1]{BCS2002}.
The maps
\begin{equation*}
	\nabla_\Gamma:H^1(\Gamma)\to\bL_t^2(\Gamma)
	\quad\text{and}\quad
	\mathbf{curl}_\Gamma:H^1(\Gamma)\to\bL_t^2(\Gamma)
\end{equation*}
define bounded and linear operators as stated in
\cite[Proposition~3.3]{BCS2002}.

Next, we set $H^{\frac32}(\Gamma)\coloneqq\gamma_\Gamma H^2(\D)$, equipped with the
quotient norm
\begin{equation*}
	\norm{\phi}_{H^{\frac32}(\Gamma)}
	\coloneqq
	\inf_{\substack{\Phi\in H^2(\D)\\ \gamma_\Gamma\Phi=\phi}}
	\norm{\Phi}_{H^2(\D)},
\end{equation*}
and denote by $H^{-\frac32}(\Gamma)$ its dual with respect to the bilinear
$L^2(\Gamma)$ duality pairing.

The operators $\nabla_\Gamma:H^{\frac32}(\Gamma)\to V_\top$ and
$\mathbf{curl}_\Gamma:H^{\frac32}(\Gamma)\to V_\times$ are bounded
\cite[Proposition~3.4]{BCS2002} and satisfy
\begin{equation*}
	\nabla_\Gamma\gamma_\Gamma\Phi=\gamma_{\top,\Gamma}\nabla\Phi
	\quad\text{and}
	\quad
	\mathbf{curl}_\Gamma\gamma_\Gamma\Phi=-\gamma_{\times,\Gamma}\nabla\Phi,
	\quad\forall \Phi\in H^2(\D).
\end{equation*}
The surface divergence $\div_\Gamma:V_\top'\to H^{-\frac32}(\Gamma)$
and the scalar surface curl
$\operatorname{curl}_\Gamma:V_\times'\to H^{-\frac32}(\Gamma)$ are the
bounded linear operators defined as (cf.~\cite[Definition~3.5]{BCS2002})
\begin{align}\label{eq:weak-surface-divergence-def}
	\dual{\div_\Gamma\bj}{\phi}_{H^{-\frac32}(\Gamma)\times H^{\frac32}(\Gamma)}
	&\coloneqq-\dual{\bj}{\nabla_\Gamma\phi}_{V_\top'\times V_\top},\\
	\label{eq:weak-surface-curl-def}
	\dual{\operatorname{curl}_\Gamma\bseta}{\phi}_{H^{-\frac32}(\Gamma)\times H^{\frac32}(\Gamma)}
	&\coloneqq\dual{\bseta}{\mathbf{curl}_\Gamma\phi}_{V_\times'\times V_\times},
\end{align}
for all $\bj\in V_\top'$, $\bseta\in V_\times'$, and $\phi\in H^{\frac32}(\Gamma)$.

Furthermore, according to \cite[Theorem~4.1]{BCS2002}, 
\begin{align}\label{eq:H-div-Gamma}
	\bH^{-\half}(\div_\Gamma,\Gamma)
	&\coloneqq\bigl\{\bj\in V_\top':
	\div_\Gamma\bj\in H^{-\half}(\Gamma)\bigr\},\\
	\label{eq:H-curl-Gamma}
	\bH^{-\half}(\operatorname{curl}_\Gamma,\Gamma)
	&\coloneqq\bigl\{\bseta\in V_\times':
	\operatorname{curl}_\Gamma\bseta\in H^{-\half}(\Gamma)\bigr\},
\end{align}
define complex Hilbert spaces equipped with the norms
\begin{equation*}
	\begin{aligned}
		\norm{\bj}_{\bH^{-\half}(\div_\Gamma,\Gamma)}^2
		&\coloneqq \norm{\bj}_{V_\top'}^2
		+\norm{\div_\Gamma\bj}_{H^{-\half}(\Gamma)}^2,\\
		\norm{\bseta}_{\bH^{-\half}(\operatorname{curl}_\Gamma,\Gamma)}^2
		&\coloneqq \norm{\bseta}_{V_\times'}^2
		+\norm{\operatorname{curl}_\Gamma\bseta}_{H^{-\half}(\Gamma)}^2.
	\end{aligned}
\end{equation*}

\subsubsection{Traces of \texorpdfstring{$\bH(\curl)$}{H(curl)} Fields and Duality}
\label{sec:hcurl-traces-duality}
The tangential trace operators $\gamma_{\times,\Gamma}$ and
$\gamma_{\top,\Gamma}$ admit unique extensions to bounded linear
surjections \cite[Theorem~4.1]{BCS2002}
\begin{equation*}
	\gamma_{\times,\Gamma}:\bH(\curl,\D)\to\bH^{-\half}(\div_\Gamma,\Gamma)
	\quad\text{and}\quad
	\gamma_{\top,\Gamma}:\bH(\curl,\D)\to\bH^{-\half}(\operatorname{curl}_\Gamma,\Gamma).
\end{equation*}
The trace operators have kernel $\bH_0(\curl,\D)$ and admit bounded
linear right inverses \cite[Section~7]{BCS2002}.

In the following, we write
\begin{equation*}
	X(\Gamma)\coloneqq\bH^{-\half}(\div_\Gamma,\Gamma)
	\quad\text{and}\quad
	Y(\Gamma)\coloneqq\bH^{-\half}(\operatorname{curl}_\Gamma,\Gamma).
\end{equation*}

We use the bilinear Maxwell pairing between $Y(\Gamma)$ and $X(\Gamma)$
with pivot space $\bL_t^2(\Gamma)$ as in \cite[Section~5, p.~858]{BCS2002}.
For $\bseta\in Y(\Gamma)\cap\bL_t^2(\Gamma)$ and
$\bj\in X(\Gamma)\cap\bL_t^2(\Gamma)$, this pairing satisfies
\begin{equation*}
	\dual{\bseta}{\bj}_{Y(\Gamma)\times X(\Gamma)}
	=\int_\Gamma\bseta\cdot \bj\,\ds.
\end{equation*}
For all $\bu,\bv\in\bH(\curl,\D)$, the Maxwell Green formula reads
\begin{equation}\label{eq:maxwell-green}
	\dual{\gamma_{\top,\Gamma}\bv}{\gamma_{\times,\Gamma}\bu}_{Y(\Gamma)\times X(\Gamma)}
	=\int_{\D}\bigl(\bv\cdot\curl\bu-\bu\cdot\curl\bv\bigr)\,\dd\x.
\end{equation}

Let $\iota_{\top,\Gamma}:V_\top'\to\bH^{-\half}(\Gamma)$ be the adjoint
operator of
$\gamma_{\top,\Gamma}\mathsf Z_\D:\bH^{\half}(\Gamma)\to V_\top$.
This bounded linear map, called the intrinsic injection, allows us to
regard a tangential density as a three-component boundary distribution
\cite[Section~2.1]{BuffaCostabelSchwab2002}. By definition,
\begin{equation*}
	\dual{\iota_{\top,\Gamma}\bj}{\bphi}_{\bH^{-\half}(\Gamma)\times\bH^{\half}(\Gamma)}
	\coloneqq\dual{\bj}{\gamma_{\top,\Gamma}\mathsf Z_\D\bphi}_{V_\top'\times V_\top},
\end{equation*}
for all $\bj\in V_\top'$ and $\bphi\in\bH^{\half}(\Gamma)$, where
$\mathsf Z_\D\bphi$ may be any $\bH^1(\D)$ lifting of $\bphi$.

The Maxwell pairing on $V_\top\times X(\Gamma)$ agrees with the
$V_\top'\times V_\top$ duality, with the arguments in reverse order.
Since $\gamma_{\top,\Gamma}\mathsf Z_\D\bphi\in V_\top$, the preceding
definition therefore gives, for all $\bj\in X(\Gamma)$ and
$\bphi\in\bH^{\half}(\Gamma)$,
\begin{equation}\label{eq:iota-top-characterization}
	\dual{\iota_{\top,\Gamma}\bj}{\bphi}_{\bH^{-\half}(\Gamma)\times
		\bH^{\half}(\Gamma)}
	=\dual{\gamma_{\top,\Gamma}\mathsf Z_{\D}\bphi}{\bj}_{Y(\Gamma)\times X(\Gamma)}.
\end{equation}

The following lemma provides a convenient characterization of the surface
divergence acting on elements of $X(\Gamma)$.

\begin{lemma}\label{lem:weak-surface-divergence}
The surface divergence operator
$\div_\Gamma\in\mathscr{L}\bigl(X(\Gamma),H^{-\half}(\Gamma)\bigr)$
satisfies, for all $\bj\in X(\Gamma)$ and $\phi\in H^{\half}(\Gamma)$,
\begin{equation}\label{eq:weak-surface-divergence}
	\dual{\div_\Gamma\bj}{\phi}_{H^{-\half}(\Gamma)\times H^{\half}(\Gamma)}
	=
	-\dual{\gamma_{\top,\Gamma}\nabla\mathsf Z_{\D}\phi}{\bj}_{Y(\Gamma)\times X(\Gamma)},
\end{equation}
where $\mathsf Z_\D\phi$ denotes an $H^1(\D)$ extension of $\phi$.
\end{lemma}

\begin{proof}
Let $\bj\in X(\Gamma)$ and $\phi\in H^{\half}(\Gamma)$. By the surjectivity
of the trace operators recalled in Subsection~\ref{sec:hcurl-traces-duality}, there exists
$\bu\in\bH(\curl,\D)$ with $\gamma_{\times,\Gamma}\bu=\bj$.

Since $\curl\bu\in\bL^2(\D)$ and $\div\curl\bu=0$,
it follows that $\curl\bu\in\bH(\div,\D)$. Its normal trace
$(\curl\bu)\cdot\bn\in H^{-\half}(\Gamma)$ therefore satisfies
\begin{equation*}
	\dual{\curl\bu\cdot\bn}{\phi}_{H^{-\half}(\Gamma)\times H^{\half}(\Gamma)}
	=
    \int_\D\curl\bu\cdot\nabla\mathsf Z_\D\phi\,\dd\x,
    \quad
    \forall 
    \phi \in H^{\frac{1}{2}}(\Gamma).
\end{equation*}
As stated in \cite[Section~4]{BCS2002},
$\div_\Gamma\bj=\div_\Gamma\gamma_{\times,\Gamma}\bu=-(\curl\bu)\cdot\bn$,
and hence
\begin{equation*}
	\dual{\div_\Gamma\bj}{\phi}_{H^{-\half}(\Gamma)\times H^{\half}(\Gamma)}
	=-\int_\D\curl\bu\cdot\nabla\mathsf Z_\D\phi\,\dd\x.
\end{equation*}

Since $\mathsf Z_\D\phi\in H^1(\D)$ and $\curl\nabla\mathsf Z_\D\phi=0$,
its gradient belongs to $\bH(\curl,\D)$, so that
$\gamma_{\top,\Gamma}\nabla\mathsf Z_\D\phi\in Y(\Gamma)$.
Applying the Maxwell Green formula \eqref{eq:maxwell-green}
with $\bv=\nabla\mathsf Z_\D\phi$ gives
\begin{align*}
	\dual{\gamma_{\top,\Gamma}\nabla\mathsf Z_\D\phi}{\bj}_{Y(\Gamma)\times X(\Gamma)}
	&=\int_\D\bigl(\curl\bu\cdot\nabla\mathsf Z_\D\phi-\bu\cdot\curl\nabla\mathsf Z_\D\phi\bigr)\,\dd\x\\
	&=\int_\D\curl\bu\cdot\nabla\mathsf Z_\D\phi\,\dd\x,
\end{align*}
thus proving \eqref{eq:weak-surface-divergence}.
\end{proof}

\subsection{Domain Transformations}
\label{sec:geometry}

Let $\widehatD\subset\IR^3$ be a fixed bounded Lipschitz domain,
in the following referred to as the \emph{reference domain},
with connected boundary $\widehatGamma\coloneqq\partial\widehatD$.
Following \cite[Subsection~2.2]{HDPS26Laplace}, we consider a collection
$\{\D_\br\}_{\br\in\frakT}$ of parametrized domains of the form
$\D_\br\coloneqq\br(\widehatD)$, for each $\br \in \mathfrak{T}$,
with boundary $\Gamma_\br\coloneqq\partial\D_\br=\br(\widehatGamma)$,
where $\frakT$ denotes the set of admissible domain transformations
and we assume that it satisfies the following conditions.
 
\begin{assumption}
\label{ass:geometry}
The set  $\frakT\subset W^{1,\infty}(\widehatD;\IR^3)$ of admissible domain transformations
satisfies the following:
\begin{itemize}
	\item[\normalfont{(i)}] 
	$\frakT$ is compact as a subset of $W^{1,\infty}(\widehatD;\IR^3)$.
	\item[\normalfont{(ii)}] 
	Each $\br\in\frakT$ is bijective and bi-Lipschitz 
	when considered as a map $\br\colon\widehatD\to\D_\br$.
	\item[\normalfont{(iii)}] 
	For each $\br\in\frakT$, the domain $\br(\widehatD)$ is a bounded Lipschitz domain in $\IR^3$.
	\item[\normalfont{(iv)}] 
	For each $\br \in \frakT$, it holds that
\begin{equation*}
	\operatorname*{ess\,inf}_{\widehat\x \in \widehatD}
	J_\br(\widehat\x)
	> 0.
\end{equation*}
\end{itemize}
where $J_\br \coloneqq \det D\br$ for each $\br \in \frakT$.
\end{assumption}

To formulate shape holomorphy, we introduce complex neighborhoods of
$\frakT$. For $\frakT$ satisfying Assumption~\ref{ass:geometry} and
$\delta>0$, we define
\begin{equation*}
	\frakT_\delta
	\coloneqq
	\left\{
	\br\in W^{1,\infty}(\widehatD;\IC^3):
	\exists\,\widetilde\br\in\frakT
	\text{ such that }
	\norm{\br-\widetilde\br}_{W^{1,\infty}(\widehatD;\IC^3)}
	<\delta
	\right\}.
\end{equation*}

\subsection{Pullback Operator}
\label{sec:pullback-operators}
Let $\mathfrak{T}$ be a set of domain transformations satisfying
Assumption~\ref{ass:geometry}.
Following \cite[Subsection~2.4]{HDPS26Laplace}, for each
$\br\in\mathfrak{T}$ we define the volume pullback operator
\begin{equation}\label{eq:volume-pullback}
	\tau_\br u\coloneqq u\circ\br,
	\quad u\in H^1(\D_\br)\ \text{or}\ u\in L^2(\D_\br).
\end{equation}
For each $\br \in \mathfrak{T}$, we consider Riesz maps
for the spaces $H^{\half}(\Gamma_\br)$ and
$H^{\half}(\widehat{\Gamma})$, i.e.,
\begin{equation}\label{eq:riesz-maps-boundary}
	\R_\br: H^{\half}(\Gamma_\br) \rightarrow H^{-\half}(\Gamma_\br)
	\quad
	\text{and}
	\quad
	\widehat{\R}: H^{\half}(\widehat{\Gamma}) \rightarrow H^{-\half}(\widehat{\Gamma}),
\end{equation}
respectively.

For convenience, we use the same notation for the boundary pullback operator,
which is defined as
\begin{equation}\label{eq:scalar-pullback}
	\tau_\br\varphi
	\coloneqq
	\left\{
	\begin{array}{cc}
		\varphi\circ\br, & \varphi\in H^{\half}(\Gamma_\br), \\
		\widehat{\R}\tau_\br\R_\br^{-1}\varphi,
		& \varphi\in H^{-\half}(\Gamma_\br),
	\end{array}
	\right.
\end{equation}
where $\R_\br$ and $\widehat{\R}$ are the Riesz maps in
\eqref{eq:riesz-maps-boundary}.
The mapping properties below follow from
\cite[Proposition~2.2]{HDPS26Laplace}.

\begin{proposition}\label{prop:pullback-isomorphisms}
Let $\mathfrak{T}$ be a set of domain transformations satisfying Assumption~\ref{ass:geometry}.
Then, for each $\br\in\mathfrak{T}$,
\begin{equation*}
	\begin{aligned}
		\tau_\br&\in
		\mathscr{L}_{\normalfont\text{iso}}(H^1(\D_\br),H^1(\widehat{\D})),
		&\quad
		\tau_\br&\in
		\mathscr{L}_{\normalfont\text{iso}}(L^2(\D_\br),L^2(\widehat{\D})),\\
		\tau_\br&\in
		\mathscr{L}_{\normalfont\text{iso}}(H^{\half}(\Gamma_\br),
		H^{\half}(\widehat{\Gamma})),
		&
		\tau_\br&\in
		\mathscr{L}_{\normalfont\text{iso}}(H^{-\half}(\Gamma_\br),
		H^{-\half}(\widehat{\Gamma})).
	\end{aligned}
\end{equation*}
In each case, the operator norms of $\tau_\br$ and $\tau_\br^{-1}$
are bounded uniformly for $\br\in\mathfrak{T}$.
\end{proposition}

For vector-valued functions, we write $\btau_\br$ for
the componentwise pullback. In the volume, it maps
$\bL^2(\D_\br)$ to $\bL^2(\widehatD)$ and
$\bH^1(\D_\br)$ to $\bH^1(\widehatD)$.
On the boundary, it maps $\bH^s(\Gamma_\br)$ to
$\bH^s(\widehatGamma)$ for $s=\pm\half$.

Next, we recall the following construction of 
the curl-conforming pullback operator from \cite{JSZ2017}.
\begin{lemma}[{\cite[Lemma~2.2]{JSZ2017}}]
\label{lem:curl-conforming-pullback}
Let $\frakT$ be a set of domain transformations satisfying
Assumption~\ref{ass:geometry}.  For each
$\br\in\frakT$, the \emph{curl-conforming pullback}
\begin{equation}\label{eq:covariant-volume-pullback}
	\curlPiola{\br}:
	\bH(\curl,\D_\br)\to\bH(\curl,\widehatD):
	\bu\mapsto(D\br)^{\mathsf T}\btau_\br\bu
\end{equation}
is a bounded isomorphism.  Its restriction to $\bH_0(\curl,\D_\br)$
is an isomorphism onto $\bH_0(\curl,\widehatD)$.
Furthermore, for all $\bu\in\bH(\curl,\D_\br)$,
\begin{equation}\label{eq:covariant-curl-transport}
	\curl(\curlPiola{\br}\bu)
	=\det(D\br)(D\br)^{-1}\btau_\br(\curl\bu)
\end{equation}
holds in $\bL^2(\widehatD)$.
\end{lemma}

The divergence-conforming counterpart of
Lemma~\ref{lem:curl-conforming-pullback} follows from the Piola identity
and a density argument, exactly as in the proof of
\cite[Lemma~2.2]{JSZ2017}.

\begin{lemma}
\label{lem:div-conforming-pullback}
Let $\frakT$ be a set of domain transformations satisfying
Assumption~\ref{ass:geometry}. For each
$\br\in\frakT$, the \emph{divergence-conforming pullback}
\begin{equation}\label{eq:covariant-volume-pullback2}
	\mathsf T_\br^{\div}:
	\bH(\div,\D_\br)\to\bH(\div,\widehatD):
	\bu\mapsto\det(D\br)(D\br)^{-1}\btau_\br\bu
\end{equation}
is a bounded isomorphism. Its restriction to $\bH_0(\div,\D_\br)$
is an isomorphism onto $\bH_0(\div,\widehatD)$.
Furthermore, for all $\bu\in\bH(\div,\D_\br)$,
\begin{equation}\label{eq:contravariant-div-transport}
	\div(\mathsf T_\br^{\div}\bu)
	=\det(D\br)\tau_\br(\div\bu)
\end{equation}
holds in $L^2(\widehatD)$.
\end{lemma}

For each $\br\in\frakT$, \eqref{eq:covariant-volume-pullback} and \eqref{eq:covariant-volume-pullback2} also define bounded isomorphisms
from $\bL^2(\D_\br)$ to $\bL^2(\widehatD)$.
For all $\bu,\bv\in\bL^2(\D_\br)$, the pullback formulas give
\begin{equation*}
	\begin{aligned}
		(\curlPiola{\br}\bu)\cdot(\mathsf T_\br^{\div}\bv)
		&=\bigl((D\br)^{\mathsf T}\btau_\br\bu\bigr)\cdot
		\bigl(\det(D\br)(D\br)^{-1}\btau_\br\bv\bigr)\\
		&=\det(D\br)\,\tau_\br(\bu\cdot\bv),
	\end{aligned}
\end{equation*}
therefore
\begin{equation}\label{eq:volume-piola-pairing}
	\int_{\widehatD}\curlPiola{\br}\bu\cdot
	\mathsf T_\br^{\div}\bv\,\dd\widehat\x
	=\int_{\D_\br}\bu\cdot\bv\,\dd\x.
\end{equation}

For $\bu\in\bL^2(\D_\br)$ and
$\widehat\bv\in\bL^2(\widehatD)$, substituting
$\bv=(\mathsf T_\br^{\div})^{-1}\widehat\bv$ in
\eqref{eq:volume-piola-pairing} gives
\begin{equation*}
	\int_{\widehatD}\curlPiola{\br}\bu\cdot\widehat\bv\,\dd\widehat\x
	=\int_{\D_\br}\bu\cdot
	(\mathsf T_\br^{\div})^{-1}\widehat\bv\,\dd\x,
\end{equation*}
therefore, for each $\br \in \mathfrak{T}$ 
\begin{equation}\label{eq:volume-piola-dual}
	(\curlPiola{\br})'=(\mathsf T_\br^{\div})^{-1}.
\end{equation}

\subsection{Covariant and Contravariant Surface Piola Transformations}
\label{sec:piola}
For the sake of convenience, for each $\br\in\frakT$,  we write $X_\br=X(\Gamma_\br)$ and
$Y_\br=Y(\Gamma_\br)$. We also set $\widehatX=X(\widehatGamma)$ and
$\widehatY=Y(\widehatGamma)$.
In the following, we construct pullback operators from $X_\br$ to
$\widehatX$ and from $Y_\br$ to $\widehatY$ to transport surface densities
from $\Gamma_\br$ to $\widehatGamma$. Following the volume-based approach
of \cite{HDPS26Helmholtz,HDPS26Laplace}, we proceed in three steps: we lift
the density in $X_\br$ to a field in $\bH(\curl,\D_\br)$, apply the curl-conforming
pullback $\curlPiola{\br}$, and, finally, take the corresponding tangential trace
on $\widehatGamma$ to obtain a surface density in $\widehat{X}$.

For $\Gamma=\partial\D$, with $\D\in\{\widehatD,\D_\br\}$ and
$\br\in\frakT$, let
\begin{equation*}
	\mathsf Z_{\D}^{\top}:Y(\Gamma)\to\bH(\curl,\D)
	\quad\text{and}\quad
	\mathsf Z_{\D}^{\times}:X(\Gamma)\to\bH(\curl,\D)
\end{equation*}
be bounded linear right inverses of the trace operators $\gamma_{\top,\Gamma}$ and $\gamma_{\times,\Gamma}$,
respectively, see \cite[Theorem~4.1]{BCS2002} for further details.
For each $\br\in\frakT$, we define the \emph{contravariant} surface Piola transformation as
\begin{equation}\label{eq:surface-piola-lifting}
	\mathsf P_\br:X_\br\to\widehatX:
	\bj\mapsto\gamma_{\times,\widehatGamma}\bigl(\curlPiola{\br}
	\mathsf Z_{\D_\br}^{\times}\bj\bigr)
\end{equation}
and the \emph{covariant} surface Piola transformation
\begin{equation}\label{eq:covariant-surface-piola-lifting}
	\mathsf Q_\br:Y_\br\to\widehatY:
	\bseta\mapsto\gamma_{\top,\widehatGamma}\bigl(\curlPiola{\br}
	\mathsf Z_{\D_\br}^{\top}\bseta\bigr).
\end{equation}

\begin{remark}
\label{rem:piola-lifting-independence}
The operators $\mathsf P_\br$ and $\mathsf Q_\br$, and the inverse
formulas stated below in Proposition~\ref{thm:piola}, are independent of the choice of
liftings. Indeed, two liftings of the same trace differ by an element of
$\bH_0(\curl,\D)$, and $\curlPiola{\br}$ and its inverse map,
as stated in Lemma~\ref{lem:curl-conforming-pullback}, the
corresponding vanishing trace spaces onto each other.
\end{remark}

\begin{proposition}
\label{thm:piola}
Let $\frakT$ be a set of domain transformations satisfying
Assumption~\ref{ass:geometry}.
For each $\br\in\frakT$, it holds that
\begin{equation}
	\mathsf P_\br\in\mathscr{L}_{\mathrm{iso}}\bigl(X_\br,\widehatX\bigr)
	\quad\text{and}\quad
	\mathsf Q_\br\in\mathscr{L}_{\mathrm{iso}}\bigl(Y_\br,\widehatY\bigr).
\end{equation}
Furthermore, for $\widehat\bj\in\widehatX$ and $\widehat{\bseta}\in\widehatY$, their inverses are given by
\begin{equation}\label{eq:surface-piola-inverses}
	\mathsf P_\br^{-1}\widehat\bj
	=\gamma_{\times,\Gamma_\br}\bigl((\curlPiola{\br})^{-1}
	\mathsf Z_{\widehatD}^{\times}\widehat\bj\bigr)
	\quad\text{and}\quad
	\mathsf Q_\br^{-1}\widehat{\bseta}
	=\gamma_{\top,\Gamma_\br}\bigl((\curlPiola{\br})^{-1}
	\mathsf Z_{\widehatD}^{\top}\widehat{\bseta}\bigr).
\end{equation}
\end{proposition}

\begin{proof}
For each $\br\in\frakT$, the lifting $\mathsf Z_{\D_\br}^{\times}:X_\br\to\bH(\curl,\D_\br)$
is bounded by construction. The map $\curlPiola{\br}:\bH(\curl,\D_\br)\to\bH(\curl,\widehatD)$
is a bounded isomorphism by Lemma~\ref{lem:curl-conforming-pullback}, and
$\gamma_{\times,\widehatGamma}:\bH(\curl,\widehatD)\to\widehatX$ is bounded.
As a composition of these three maps, we have that $\mathsf P_\br\in\mathscr L(X_\br,\widehatX)$. The same argument
with $\mathsf Z_{\D_\br}^{\top}$ and $\gamma_{\top,\widehatGamma}$ gives
$\mathsf Q_\br\in\mathscr L(Y_\br,\widehatY)$. 

Similarly, each operator in \eqref{eq:surface-piola-inverses} is the 
composition of a lifting on $\widehatD$, the isomorphism $(\curlPiola{\br})^{-1}$, and a
trace on $\Gamma_\br$. Hence they define bounded linear maps from
$\widehatX$ to $X_\br$ and from $\widehatY$ to $Y_\br$, respectively.

We next verify that the map in \eqref{eq:surface-piola-inverses} is a
left inverse of $\mathsf P_\br$. Let $\bj\in X_\br$ and set
$\bu\coloneqq\mathsf Z_{\D_\br}^{\times}\bj\in\bH(\curl,\D_\br)$, so that
$\gamma_{\times,\Gamma_\br}\bu=\bj$ and, by
\eqref{eq:surface-piola-lifting},
$\mathsf P_\br\bj=\gamma_{\times,\widehatGamma}(\curlPiola{\br}\bu)$.
The field $\curlPiola{\br}\bu\in\bH(\curl,\widehatD)$ is therefore a
lifting of $\mathsf P_\br\bj$, and using it in place of
$\mathsf Z_{\widehatD}^{\times}\mathsf P_\br\bj$ in
\eqref{eq:surface-piola-inverses} (cf. Remark~\ref{rem:piola-lifting-independence}) yields
\begin{equation*}
	\mathsf P_\br^{-1}\mathsf P_\br\bj
	=\gamma_{\times,\Gamma_\br}\bigl((\curlPiola{\br})^{-1}\curlPiola{\br}\bu\bigr)
	=\gamma_{\times,\Gamma_\br}\bu
	=\bj.
\end{equation*}

Let $\widehat\bj\in\widehatX$ and set
$\widehat\bu\coloneqq\mathsf Z_{\widehatD}^{\times}\widehat\bj\in\bH(\curl,\widehatD)$,
so that $\gamma_{\times,\widehatGamma}\widehat\bu=\widehat\bj$ and
\eqref{eq:surface-piola-inverses} reads
$\mathsf P_\br^{-1}\widehat\bj=\gamma_{\times,\Gamma_\br}((\curlPiola{\br})^{-1}\widehat\bu)$.
Now $(\curlPiola{\br})^{-1}\widehat\bu\in\bH(\curl,\D_\br)$ is a lifting
of $\mathsf P_\br^{-1}\widehat\bj$, and using it in
\eqref{eq:surface-piola-lifting} (cf. Remark~\ref{rem:piola-lifting-independence}) yields
\begin{equation*}
	\mathsf P_\br\mathsf P_\br^{-1}\widehat\bj
	=\gamma_{\times,\widehatGamma}\bigl(\curlPiola{\br}(\curlPiola{\br})^{-1}\widehat\bu\bigr)
	=\gamma_{\times,\widehatGamma}\widehat\bu
	=\widehat\bj.
\end{equation*}
Consequently, $\mathsf P_\br$ is a boundedly invertible isomorphism with inverse
given by \eqref{eq:surface-piola-inverses}, that is,
$\mathsf P_\br\in\mathscr{L}_{\mathrm{iso}}(X_\br,\widehatX)$.

The argument for $\mathsf Q_\br$ follows along the same lines, with
$\gamma_{\times,\Gamma}$ and $\mathsf Z^{\times}_{\D}$ replaced by
$\gamma_{\top,\Gamma}$ and $\mathsf Z^{\top}_{\D}$ on each of the domains
$\D\in\{\D_\br,\widehatD\}$, where $\Gamma=\partial\D$.
In both compositions, the transported field serves as the lifting in the
second map (cf. Remark~\ref{rem:piola-lifting-independence}) and the
original trace is recovered.
\end{proof}

\begin{lemma}
\label{lem:piola-duality}
Let $\frakT$ be a set of domain transformations satisfying
Assumption~\ref{ass:geometry}.
For each $\br\in\frakT$, for all $\bj\in X_\br$ and
 $\bseta\in Y_\br$, it holds that
\begin{equation}\label{eq:intrinsic-surface-piola}
	\dual{\mathsf Q_\br\bseta}{\mathsf P_\br\bj}_{\widehatY\times\widehatX}
	=\dual{\bseta}{\bj}_{Y_\br\times X_\br}.
\end{equation}
\end{lemma}

\begin{proof}
For each $\br\in\frakT$, let $\bj\in X_\br$ and $\bseta\in Y_\br$.
Set $\bu=\mathsf Z_{\D_\br}^{\times}\bj$,
$\bv=\mathsf Z_{\D_\br}^{\top}\bseta$, as well as $\widehat\bu=\curlPiola{\br}\bu$
and $\widehat\bv=\curlPiola{\br}\bv$.
By definition, $\mathsf P_\br\bj=\gamma_{\times,\widehatGamma} \widehat \bu$
and $\mathsf Q_\br\bseta=\gamma_{\top,\widehatGamma}\widehat\bv$, and
it follows from  \eqref{eq:covariant-curl-transport} that
$\curl \widehat\bu=\mathsf T_\br^{\div}(\curl\bu)$ and
$\curl\widehat\bv=\mathsf T_\br^{\div}(\curl\bv)$.

The application of \eqref{eq:volume-piola-pairing} together with the Maxwell Green formula
stated in \eqref{eq:maxwell-green} yields
\begin{align*}
	\dual{
		\mathsf Q_\br\bseta
	}{
		\mathsf P_\br\bj
	}_{\widehatY\times\widehatX}
	&
	=
	\int_{\widehatD}
	\left(
	\widehat\bv\cdot\curl \widehat \bu
	-
	\widehat \bu \cdot\curl\widehat\bv
	\right)\,\dd\widehat\x\\
	&
	=
	\int_{\widehatD}
	\bigl(
	(\curlPiola{\br}\bv)\cdot\mathsf T_\br^{\div}(\curl\bu)
	-(\curlPiola{\br}\bu)\cdot\mathsf T_\br^{\div}(\curl\bv)
	\bigr)\,\dd\widehat\x\\
	&
	=
	\int_{\D_\br}\bigl(\bv\cdot\curl\bu-\bu\cdot\curl\bv\bigr)\,\dd\x\\
	&
	=
	\dual{\bseta}{\bj}_{Y_\br\times X_\br}.
\end{align*}
\end{proof}

The next result addresses how the surface divergence transforms under
the pullback operator $\mathsf P_\br$. 

\begin{lemma}
\label{lem:piola-div}
Let $\frakT$ be a set of domain transformations satisfying
Assumption~\ref{ass:geometry}.
For each $\br\in\frakT$ and for all $\bj\in X_\br$ and
 $\phi\in H^{\half}(\Gamma_\br)$, it holds that
\begin{equation}\label{eq:piola-div}
	\dual{
		\div_{\Gamma_\br}\bj
	}{
		\phi
	}_{
		H^{-\half}(\Gamma_\br)\times H^{\half}(\Gamma_\br)
	}
	=
	\dual{
		\divGh(\mathsf P_\br\bj)
	}{
		\tau_\br\phi
	}_{
		H^{-\half}(\widehatGamma)\times H^{\half}(\widehatGamma)
	}.
\end{equation}
\end{lemma}

\begin{proof}
For each $\br\in\frakT$, consider
$\bj\in X_\br$ and $\phi\in H^{\half}(\Gamma_\br)$, and set
$\widehat\bj=\mathsf P_\br\bj$ and $\widehat\phi=\tau_\br\phi$.
According to Lemma~\ref{lem:weak-surface-divergence}, for each $\br \in \mathfrak{T}$, we
have that
\begin{equation*}
	\dual{\div_{\Gamma_\br}\bj}{\phi}_{H^{-\half}(\Gamma_\br)\times H^{\half}(\Gamma_\br)}
	=
	-\dual{\gamma_{\top,\Gamma_\br}\nabla\mathsf Z_{\D_\br}\phi}{\bj}_{Y_\br\times X_\br}.
\end{equation*}
Observe that, by construction, $\tau_\br\mathsf Z_{\D_\br}\phi\in H^1(\widehatD)$ has
Dirichlet trace $\widehat\phi$.
Recalling Remark~\ref{rem:piola-lifting-independence}, we may use $\nabla\mathsf Z_{\D_\br}\phi$
as the lifting of $\gamma_{\top,\Gamma_\br}\nabla\mathsf Z_{\D_\br}\phi$ in the definition
of $\mathsf Q_\br$. By \eqref{eq:covariant-volume-pullback} and the
Sobolev chain rule in the volume, we have
\begin{equation*}
	\curlPiola{\br}\nabla\mathsf Z_{\D_\br}\phi
	=(D\br)^{\mathsf T}\btau_\br(\nabla\mathsf Z_{\D_\br}\phi)
	=\widehat{\nabla}(\tau_\br\mathsf Z_{\D_\br}\phi).
\end{equation*}
Taking the tangential trace on $\widehatGamma$ gives
\begin{equation*}
	\mathsf Q_\br(\gamma_{\top,\Gamma_\br}\nabla\mathsf Z_{\D_\br}\phi)
	=\gamma_{\top,\widehatGamma}(\curlPiola{\br}\nabla\mathsf Z_{\D_\br}\phi)
	=\gamma_{\top,\widehatGamma}\widehat{\nabla}(\tau_\br\mathsf Z_{\D_\br}\phi).
\end{equation*}
Applying \eqref{eq:intrinsic-surface-piola} with
$\bseta=\gamma_{\top,\Gamma_\br}\nabla\mathsf Z_{\D_\br}\phi$ yields
\begin{align*}
	\dual{\div_{\Gamma_\br}\bj}{\phi}_{H^{-\half}(\Gamma_\br)\times H^{\half}(\Gamma_\br)}
	&=-\dual{\gamma_{\top,\widehatGamma}\widehat{\nabla}(\tau_\br\mathsf Z_{\D_\br}\phi)}{\widehat\bj}_{\widehatY\times\widehatX}\\
	&=\dual{\divGh\widehat\bj}{\widehat\phi}_{H^{-\half}(\widehatGamma)\times H^{\half}(\widehatGamma)},
\end{align*}
where the last equality follows from Lemma~\ref{lem:weak-surface-divergence}
on $\widehatGamma$ with the lifting $\tau_\br\mathsf Z_{\D_\br}\phi$.
\end{proof}

\begin{remark}
\label{rem:classical-surface-piola}
Suppose that the shape transformation $\br$, the boundaries, and the tangential fields $\bj,\bseta$
on $\Gamma_\br$ are smooth. From the definitions
\eqref{eq:surface-piola-lifting} and \eqref{eq:covariant-surface-piola-lifting},
we obtain the classical surface Piola formulas
\begingroup\postdisplaypenalty=10000\relax
\begin{equation}\label{eq:classical-contravariant-surface-piola}
	(\mathsf P_\br\bj)(\widehat{\x})
	=J_{\mathrm S,\br}(\widehat{\x})
	\bigl(d_{\widehatGamma}\br(\widehat{\x})\bigr)^{-1}
	(\bj\circ\br)(\widehat{\x}),\quad \widehat{\x}\in\widehatGamma,
\end{equation}
and
\begin{equation}\label{eq:classical-covariant-surface-piola}
	(\mathsf Q_\br\bseta)(\widehat{\x})
	=\bigl(d_{\widehatGamma}\br(\widehat{\x})\bigr)^{\mathsf T}
	(\bseta\circ\br)(\widehat{\x}),\quad \widehat{\x}\in\widehatGamma,
\end{equation}
\endgroup
where $d_{\widehatGamma}\br(\widehat{\x})$ denotes the restriction
$\widehat\bv\mapsto D\br(\widehat{\x})\widehat\bv$ to vectors satisfying
$\widehat\bv\cdot\widehat{\bn}(\widehat{\x})=0$,
which maps the plane orthogonal to $\widehat{\bn}(\widehat{\x})$
onto the plane orthogonal to $\bn(\br(\widehat{\x}))$, and
$J_{\mathrm S,\br}$ is the surface Jacobian of $\br|_{\widehatGamma}$.

Let $\bv$ be a smooth lifting of $\bj$ such that
$\gamma_{\times,\Gamma_\br}\bv=\bj$.
At $\widehat{\x}$, choose an orthonormal
tangent basis $\mathbf e_1,\mathbf e_2$ with
$\mathbf e_1\times\mathbf e_2=\widehat{\bn}$. Recalling the identity
$(A\mathbf a)\times(A\mathbf b)=\operatorname{cof}(A)(\mathbf a\times\mathbf b)$,
we obtain
\begin{equation*}
	(D\br\,\mathbf e_1)\times(D\br\,\mathbf e_2)
	=\operatorname{cof}(D\br)\widehat{\bn}
	=J_{\mathrm S,\br}\btau_\br\bn .
\end{equation*}
Recalling also
$A(\mathbf{a}\times A^{\mathsf T}\mathbf{b})=(\operatorname{cof}(A)\mathbf{a})\times\mathbf{b}$,
we apply $D\br$ to \eqref{eq:surface-piola-lifting}, using $\bv$ as the lifting
(cf. Remark~\ref{rem:piola-lifting-independence}), and obtain
\begin{equation*}
	d_{\widehatGamma}\br\,\mathsf P_\br\bj
	=D\br\bigl(\widehat{\bn}\times (D\br)^{\mathsf T}\btau_\br\bv\bigr)
	=(\operatorname{cof}(D\br)\widehat{\bn})\times\btau_\br\bv.
\end{equation*}
Substituting $\operatorname{cof}(D\br)\widehat{\bn}
=J_{\mathrm S,\br}\btau_\br\bn$ gives
\begin{equation*}
	(\operatorname{cof}(D\br)\widehat{\bn})\times\btau_\br\bv
	=J_{\mathrm S,\br}\btau_\br(\bn\times\bv)
	=J_{\mathrm S,\br}\btau_\br\bj .
\end{equation*}
Applying $(d_{\widehatGamma}\br)^{-1}$ gives the formula for $\mathsf P_\br$
in \eqref{eq:classical-contravariant-surface-piola}.

Let $\bu$ be a smooth lifting of $\bseta$, with
$\gamma_{\top,\Gamma_\br}\bu=\bseta$. By
Remark~\ref{rem:piola-lifting-independence}, we may use $\bu$ in
\eqref{eq:covariant-surface-piola-lifting}, so
\eqref{eq:covariant-volume-pullback} identifies $\mathsf Q_\br\bseta$ as the
tangential projection of $(D\br)^{\mathsf T}\btau_\br\bu$ on
$\widehatGamma$, namely,
\begin{equation*}
	\mathsf Q_\br\bseta
	=(D\br)^{\mathsf T}\btau_\br\bu
	-\bigl(\bigl((D\br)^{\mathsf T}\btau_\br\bu\bigr)\cdot\widehat{\bn}\bigr)\widehat{\bn}.
\end{equation*}
Taking the dot product with a tangent vector $\widehat\bv$ eliminates the
normal term because $\widehat{\bn}\cdot\widehat\bv=0$.
Moving $(D\br)^{\mathsf T}$ across the dot product and using
$D\br\,\widehat\bv=d_{\widehatGamma}\br\,\widehat\bv$ then gives
\begin{equation*}
	(\mathsf Q_\br\bseta)\cdot\widehat\bv
	=\bigl((D\br)^{\mathsf T}\btau_\br\bu\bigr)\cdot\widehat\bv
	=(\btau_\br\bu)\cdot(d_{\widehatGamma}\br\,\widehat\bv).
\end{equation*}
Since $\br$ maps $\widehatGamma$ onto $\Gamma_\br$, the vector
$d_{\widehatGamma}\br\,\widehat\bv$ is tangent to $\Gamma_\br$
at $\br(\widehat{\x})$. This dot product therefore depends only on the
tangential component $\bseta$ of $\bu$, giving
\begin{equation*}
	(\btau_\br\bu)\cdot(d_{\widehatGamma}\br\,\widehat\bv)
	=(\btau_\br\bseta)\cdot(d_{\widehatGamma}\br\,\widehat\bv)
	=\bigl((d_{\widehatGamma}\br)^{\mathsf T}\btau_\br\bseta\bigr)
	\cdot\widehat\bv.
\end{equation*}
Both $\mathsf Q_\br\bseta$ and
$(d_{\widehatGamma}\br)^{\mathsf T}\btau_\br\bseta$ are tangential to
$\widehatGamma$. Since their dot products agree for every tangent vector
$\widehat\bv$, they coincide, proving
\eqref{eq:classical-covariant-surface-piola}.
\end{remark}

\section{Weak Transport to the Reference Boundary}
\label{sec:transport}

To express the EFIE on the reference boundary, we need the
transported intrinsic injection and surface divergence.
We establish their holomorphic dependence on the deformation
through the transported duality pairings.

For each $\br\in\frakT$, we introduce the bilinear form
$\dual{\cdot}{\cdot}_{H^{-\half}(\widehatGamma)\times H^{\half}(\widehatGamma),\br}:
 H^{-\half}(\widehatGamma)\times H^{\half}(\widehatGamma)\to\IC$
as in \cite[Subsection~2.4]{HDPS26Laplace}.
For all $\widehat \varphi\in H^{-\half}(\widehatGamma)$ and
$\widehat\phi\in H^{\half}(\widehatGamma)$, it is given by
\begin{equation}\label{eq:definition-transported-duality}
	\begin{aligned}
		&\dual{\widehat \varphi}{\widehat\phi}_{H^{-\half}(\widehatGamma)\times H^{\half}(\widehatGamma),\br}
		\\
		&\quad\coloneqq
		\int_{\widehatGamma}
		(\widehat{\R}^{-1}\widehat \varphi)\widehat\phi
		J_{\mathrm S,\br}\ds_{\widehat{\x}}
		\\
		&\qquad{}+
		\int_{\widehatGamma}\int_{\widehatGamma}
		\frac{\bigl((\widehat{\R}^{-1}\widehat \varphi)
			(\widehat{\x})
			-(\widehat{\R}^{-1}\widehat \varphi)
			(\widehat{\y})\bigr)
			\bigl(\widehat\phi(\widehat{\x})
			-\widehat\phi(\widehat{\y})\bigr)}
		{\norm{\br(\widehat{\x})
				-\br(\widehat{\y})}^3}
		J_{\mathrm S,\br}(\widehat{\x})
		J_{\mathrm S,\br}(\widehat{\y})
		\ds_{\widehat{\x}}\ds_{\widehat{\y}}.
	\end{aligned}
\end{equation}
Furthermore, for each $\br\in\frakT$ and for all
$\varphi\in H^{-\half}(\Gamma_\br)$, $\phi\in H^{\half}(\Gamma_\br)$,
according to \cite[Subsection~2.4]{HDPS26Laplace}, we have that
\begin{equation}\label{eq:transported-scalar-duality}
	\dual{\varphi}{\phi}_{H^{-\half}(\Gamma_\br)\times H^{\half}(\Gamma_\br)}
	=\dual{\tau_\br \varphi}{\tau_\br\phi}_{H^{-\half}(\widehatGamma)\times H^{\half}(\widehatGamma),\br}.
\end{equation}
We recall the following result from
\cite[Lemma~2.7]{HDPS26Laplace}.
\begin{lemma}
\label{lem:transported-duality}
Let $\mathfrak{T}$ be a set of domain transformations
satisfying Assumption \ref{ass:geometry}.
\begin{itemize}
\item[\normalfont{(i)}]
There exists $\delta>0$ such that 
the bilinear form 
$ \dual{\cdot}{\cdot}_{{H}^{-\frac{1}{2}}(\widehat\Gamma) \times H^\frac{1}{2}(\widehat\Gamma),\br}:
{H}^{-\frac{1}{2}}(\widehat\Gamma) \times H^\frac{1}{2}(\widehat\Gamma) \rightarrow \mathbb{C}$ 
introduced in \eqref{eq:definition-transported-duality}
admits a bounded extension to $\mathfrak{T}_\delta$ .
More precisely, for each $\br \in \mathfrak{T}_\delta$ 
and for all $\widehat{\phi} \in {H}^{-\frac{1}{2}}(\widehat\Gamma)$,
$\widehat{\varphi}  \in H^{\frac{1}{2}}(\widehat\Gamma)$, it holds that
\begin{equation}\label{eq:transported-duality-bound}
	\snorm{
		\dual{
			\widehat{\phi}
		}{
			\widehat{\varphi}
		}_{{H}^{-\frac{1}{2}}(\widehat\Gamma) \times H^{\frac{1}{2}}(\widehat\Gamma),\br}
	}
	\leq
	C_1(\mathfrak{T},\delta)
	\norm{\widehat{\phi} }_{{H}^{-\frac{1}{2}}(\widehat\Gamma)}
	\norm{\widehat{\varphi} }_{{H}^{\frac{1}{2}}(\widehat\Gamma)},
\end{equation}
where $C_1(\mathfrak{T},\delta)>0$ depends only on $\mathfrak{T}$
and $\delta>0$.
\item[\normalfont{(ii)}]
There exists $\delta>0$ such that the map
\begin{equation}\label{eq:transported-duality-holomorphy}
	\mathfrak{T} \rightarrow
	\mathscr{L}_{\normalfont\text{bl}}
	\left(
	H^{-\frac{1}{2}}(\widehat\Gamma)  \times H^{\frac{1}{2}}(\widehat\Gamma)
	,
	\mathbb{C}
	\right):
	\br
	\mapsto
	\left(
	(\widehat\phi,\widehat\psi)
	\mapsto
	\dual{\widehat\phi}{\widehat\psi}_{H^{-\frac{1}{2}}(\widehat\Gamma)  \times H^{\frac{1}{2}}(\widehat\Gamma),\br}
	\right)
\end{equation}
admits a bounded holomorphic extension to $\mathfrak{T}_\delta$.
\item[\normalfont{(iii)}]
There exist $\delta>0$ and $C_2(\mathfrak{T},\delta)>0$, depending only
on $\mathfrak{T}$ and $\delta$, such that
for each $\br \in \mathfrak{T}_\delta$
\begin{equation}\label{eq:transported-duality-infsup}
	\begin{aligned}
		\inf_{\widehat{\phi}\in {H}^{-\frac{1}{2}}(\widehat{\Gamma})\setminus\{0\}}
		\sup_{\widehat{\psi}\in H^{\frac{1}{2}}(\widehat{\Gamma})\setminus\{0\}}
		\frac{
			\snorm{\dual{\widehat\phi}{\widehat\psi}_{H^{-\frac{1}{2}}(\widehat\Gamma)  \times H^{\frac{1}{2}}(\widehat\Gamma),\br}}
		}{
			\norm{\widehat{\phi}}_{{H}^{-\frac{1}{2}}(\widehat{\Gamma})}
			\norm{\widehat{\psi}}_{{H}^{\frac{1}{2}}(\widehat{\Gamma})}
		}
		&
		\geq C_2(\mathfrak{T},\delta)
		\quad
		\text{and}
		\\
		\inf_{\widehat{\psi}\in H^{\frac{1}{2}}(\widehat{\Gamma})\setminus\{0\}}
		\sup_{\widehat{\phi}\in {H}^{-\frac{1}{2}}(\widehat{\Gamma})\setminus\{0\}}
		\frac{
			\snorm{\dual{\widehat\phi}{\widehat\psi}_{H^{-\frac{1}{2}}(\widehat\Gamma)  \times H^{\frac{1}{2}}(\widehat\Gamma),\br}}
		}{
			\norm{\widehat{\phi}}_{{H}^{-\frac{1}{2}}(\widehat{\Gamma})}
			\norm{\widehat{\psi}}_{{H}^{\frac{1}{2}}(\widehat{\Gamma})}
		}
		&
		\geq C_2(\mathfrak{T},\delta).
	\end{aligned}
\end{equation}
\end{itemize}
\end{lemma}

For each $\br\in\frakT$ and for all
$\widehat{\boldsymbol\varphi}\in\bH^{-\half}(\widehatGamma)$ and
$\widehat{\bphi}\in\bH^{\half}(\widehatGamma)$, we define
\begin{equation}\label{eq:transported-scalar-duality-vector}
	\dual{
		\widehat{\boldsymbol\varphi}
	}{
		\widehat{\bphi}
	}_ {\bH^{-\half}(\widehatGamma)\times\bH^{\half}(\widehatGamma),\br}
	\coloneqq
	\sum_{\ell=1}^{3}
	\dual{\widehat\varphi_\ell}{\widehat\phi_\ell}_
	{H^{-\half}(\widehatGamma)\times H^{\half}(\widehatGamma),\br}.
\end{equation}

The next result follows componentwise from
Lemma~\ref{lem:transported-duality}.
\begin{corollary}
\label{cor:transported-vector-duality}
Let $\frakT$ be a set of domain transformations satisfying
Assumption~\ref{ass:geometry}.
\begin{itemize}
\item[\normalfont{(i)}]
There exists $\delta>0$ such that the bilinear form
$\dual{\cdot}{\cdot}_{\bH^{-\half}(\widehatGamma)\times\bH^{\half}(\widehatGamma),\br}:
\bH^{-\half}(\widehatGamma)\times\bH^{\half}(\widehatGamma)\to\IC$
introduced above admits a bounded extension to $\frakT_\delta$.
More precisely, for each $\br\in\frakT_\delta$ and for all
$\widehat{\boldsymbol\varphi}\in\bH^{-\half}(\widehatGamma)$,
$\widehat{\bphi}\in\bH^{\half}(\widehatGamma)$, it holds that
\begin{equation*}
	\snorm{\dual{\widehat{\boldsymbol\varphi}}{\widehat{\bphi}}_
		{\bH^{-\half}(\widehatGamma)\times\bH^{\half}(\widehatGamma),\br}}
	\leq C_1(\frakT,\delta)
	\norm{\widehat{\boldsymbol\varphi}}_{\bH^{-\half}(\widehatGamma)}
	\norm{\widehat{\bphi}}_{\bH^{\half}(\widehatGamma)},
\end{equation*}
where $C_1(\frakT,\delta)>0$ depends only on $\frakT$
and $\delta>0$ and is as in \cref{eq:transported-duality-bound}.
\item[\normalfont{(ii)}]
There exists $\delta>0$ such that the map
\begin{equation*}
	\frakT\to\mathscr{L}_{\normalfont\text{bl}}
	\bigl(\bH^{-\half}(\widehatGamma)\times\bH^{\half}(\widehatGamma),\IC\bigr):
	\br\mapsto\left(
	(\widehat{\boldsymbol\varphi},\widehat{\bphi})\mapsto
	\dual{\widehat{\boldsymbol\varphi}}{\widehat{\bphi}}_
	{\bH^{-\half}(\widehatGamma)\times\bH^{\half}(\widehatGamma),\br}
	\right)
\end{equation*}
admits a bounded holomorphic extension to $\frakT_\delta$.
\item[\normalfont{(iii)}]
There exist $\delta>0$ and $C_2(\frakT,\delta)>0$, depending only
on $\frakT$ and $\delta$, such that
for each $\br\in\frakT_\delta$
\begin{equation*}
	\begin{aligned}
		\inf_{\widehat{\boldsymbol\varphi}\in\bH^{-\half}(\widehatGamma)\setminus\{0\}}
		\sup_{\widehat{\bphi}\in\bH^{\half}(\widehatGamma)\setminus\{0\}}
		\frac{\snorm{\dual{\widehat{\boldsymbol\varphi}}{\widehat{\bphi}}_
				{\bH^{-\half}(\widehatGamma)\times\bH^{\half}(\widehatGamma),\br}}}
		{\norm{\widehat{\boldsymbol\varphi}}_{\bH^{-\half}(\widehatGamma)}
			\norm{\widehat{\bphi}}_{\bH^{\half}(\widehatGamma)}}
		&\geq C_2(\frakT,\delta)
		\quad\text{and}\\
		\inf_{\widehat{\bphi}\in\bH^{\half}(\widehatGamma)\setminus\{0\}}
		\sup_{\widehat{\boldsymbol\varphi}\in\bH^{-\half}(\widehatGamma)\setminus\{0\}}
		\frac{\snorm{\dual{\widehat{\boldsymbol\varphi}}{\widehat{\bphi}}_
				{\bH^{-\half}(\widehatGamma)\times\bH^{\half}(\widehatGamma),\br}}}
		{\norm{\widehat{\boldsymbol\varphi}}_{\bH^{-\half}(\widehatGamma)}
			\norm{\widehat{\bphi}}_{\bH^{\half}(\widehatGamma)}}
		&\geq C_2(\frakT,\delta),
	\end{aligned}
\end{equation*}
with $C_2(\frakT,\delta)$ as in Lemma~\ref{lem:transported-duality}.
\end{itemize}
\end{corollary}

The next result is the vector-valued extension of
\cite[Lemma~2.8]{HDPS26Laplace}.
\begin{lemma}
\label{lem:vector-variational-solution}
Let $\frakT$ be a set of domain transformations
satisfying Assumption~\ref{ass:geometry}.
Assume that there exists $\varepsilon>0$ such that the map
$\frakT\ni\br\mapsto\bff_\br\in\bH^{-\half}(\widehatGamma)$
admits a bounded holomorphic extension to $\frakT_\varepsilon$.
Then the following assertions hold.
\begin{itemize}
\item[\normalfont{(i)}]
There exists $\delta\in(0,\varepsilon)$ such that for each $\br\in\frakT_\delta$
there exists a unique $\widehat{\bpsi}_\br\in\bH^{-\half}(\widehatGamma)$
satisfying
\begin{equation}\label{eq:vector-variational-problem}
	\dual{\widehat{\bpsi}_\br}{\widehat{\bphi}}_
	{\bH^{-\half}(\widehatGamma)\times\bH^{\half}(\widehatGamma),\br}
	=\bff_\br(\widehat{\bphi}),
	\quad\forall\widehat{\bphi}\in\bH^{\half}(\widehatGamma).
\end{equation}
\item[\normalfont{(ii)}]
There exists $\delta\in(0,\varepsilon)$ such that the map
\begin{equation*}
	\frakT\to\bH^{-\half}(\widehatGamma):
	\br\mapsto\widehat{\bpsi}_\br
\end{equation*}
admits a bounded holomorphic extension to $\frakT_\delta$.
\end{itemize}
\end{lemma}

For use in the subsequent analysis of the EFIE, we define, for each
$\br\in\frakT$, the transported intrinsic injection and surface divergence as
\begin{equation}\label{eq:coefficient-identification}
	\widehat{\iota}_{\top,\br}
	\coloneqq\btau_\br\iota_{\top,\Gamma_\br}\mathsf P_\br^{-1}
	\in\mathscr L\bigl(\widehatX,\bH^{-\half}(\widehatGamma)\bigr)
\end{equation}
and
\begin{equation}\label{eq:transported-surface-divergence}
	\widehat{\div}_\br
	\coloneqq\tau_\br\div_{\Gamma_\br}\mathsf P_\br^{-1}
	\in\mathscr L\bigl(\widehatX,H^{-\half}(\widehatGamma)\bigr).
\end{equation}
We next establish the holomorphic dependence of these two operators
on the domain transformation.

\begin{proposition}
\label{prop:coefficient-identification}
Let $\frakT$ be a set of domain transformations satisfying
Assumption~\ref{ass:geometry}. Then there exists $\delta>0$ such that the maps
\begin{equation}\label{eq:coefficient-holomorphy}
	\mathcal I:\frakT\to
	\mathscr{L}\bigl(\widehatX,\bH^{-\half}(\widehatGamma)\bigr):
	\br\mapsto\widehat{\iota}_{\top,\br}
\end{equation}
and
\begin{equation}\label{eq:div-coefficient-holomorphy}
	\mathcal D:\frakT\to
	\mathscr{L}\bigl(\widehatX,H^{-\half}(\widehatGamma)\bigr):
	\br\mapsto\widehat{\div}_\br
\end{equation}
admit bounded holomorphic extensions to $\frakT_\delta$.
\end{proposition}

\begin{proof}
Let $\widehat\bj\in\widehatX$ and
$\widehat{\bphi}\in\bH^{\half}(\widehatGamma)$,
and for each $\br\in\frakT$ set
$\bj_\br=\mathsf P_\br^{-1}\widehat\bj\in X_\br$ and
$\bphi_\br=\btau_\br^{-1}\widehat{\bphi}
\in\bH^{\half}(\Gamma_\br)$ together with
$\bpsi_\br=\btau_\br^{-1}\mathsf Z_{\widehatD}\widehat{\bphi}
\in\bH^1(\D_\br)$, where $\mathsf Z_{\widehatD}$ is the componentwise
Dirichlet lifting.
The Dirichlet trace commutes with the pullback operator, thus
\begin{equation*}
	\gamma_{\Gamma_\br}\bpsi_\br
	=
	\btau_\br^{-1}\gamma_{\widehatGamma}
	\mathsf Z_{\widehatD}\widehat{\bphi}
	=
	\bphi_\br.
\end{equation*}
Using \eqref{eq:iota-top-characterization},
\eqref{eq:transported-scalar-duality},
\eqref{eq:transported-scalar-duality-vector},
and \eqref{eq:coefficient-identification}, we obtain
\begin{align}\label{eq:transport_injection}
	\dual{
		\widehat{\iota}_{\top,\br}\widehat\bj
	}{
		\widehat{\bphi}
	}_{\bH^{-\half}(\widehatGamma)\times\bH^{\half}(\widehatGamma),\br}
	&=
	\dual{
		\iota_{\top,\Gamma_\br}\bj_\br
	}{
		\bphi_\br
	}_{\bH^{-\half}(\Gamma_\br)\times\bH^{\half}(\Gamma_\br)}\\
	&=
	\dual{\gamma_{\top,\Gamma_\br}\bpsi_\br}
	{\bj_\br}_{Y_\br\times X_\br},
\end{align}
where we have used $\bpsi_\br$ as the lifting of $\bphi_\br$.
Next, observe that
\begin{equation*}
	\curlPiola{\br}\bpsi_\br
	=(D\br)^{\mathsf T}\btau_\br\bpsi_\br
	=(D\br)^{\mathsf T}\mathsf Z_{\widehatD}\widehat{\bphi}
	\quad\text{in }\widehatD,
\end{equation*}
where we have used that by construction
$\btau_\br\bpsi_\br=\mathsf Z_{\widehatD}\widehat{\bphi}$.
Recalling Remark~\ref{rem:piola-lifting-independence}, we therefore use
$\bpsi_\br\in\bH^1(\D_\br)\subset\bH(\curl,\D_\br)$ as the lifting of
$\gamma_{\top,\Gamma_\br}\bpsi_\br$ in the definition of $\mathsf Q_\br$.
Lemma~\ref{lem:piola-duality} and \eqref{eq:transport_injection} yield
\begin{equation}\label{eq:coefficient-reference-pairing}
	\begin{aligned}
		\dual{\widehat{\iota}_{\top,\br}\widehat\bj}{\widehat{\bphi}}_
		{\bH^{-\half}(\widehatGamma)\times\bH^{\half}(\widehatGamma),\br}
		&=
		\dual{\mathsf Q_\br(\gamma_{\top,\Gamma_\br}\bpsi_\br)}
		{\mathsf P_\br\bj_\br}_{\widehatY\times\widehatX}\\
		&=
		\dual{\gamma_{\top,\widehatGamma}(\curlPiola{\br}\bpsi_\br)}
		{\widehat\bj}_{\widehatY\times\widehatX}\\
		&=
		\dual{\gamma_{\top,\widehatGamma}\!\left((D\br)^{\mathsf T}
			\mathsf Z_{\widehatD}\widehat{\bphi}\right)}
		{\widehat\bj}_{\widehatY\times\widehatX}.
	\end{aligned}
\end{equation}

For $\br\in W^{1,\infty}(\widehatD;\IC^3)$ and
$\widehat\bw\in\bH^1(\widehatD)$, let $r_\ell$ and $\widehat w_\ell$
denote the components of $\br$ and $\widehat\bw$, respectively.
By the product rule and an approximation argument as in the proof of
\cite[Lemma~2.2]{JSZ2017}, using $\curl\widehat{\nabla} r_\ell=0$, we obtain
\begin{equation*}
	\curl\bigl((D\br)^{\mathsf T}\widehat\bw\bigr)
	=\curl\Bigl(\sum_{\ell=1}^3\widehat w_\ell\widehat{\nabla} r_\ell\Bigr)
	=\sum_{\ell=1}^3\widehat{\nabla}\widehat w_\ell\times\widehat{\nabla} r_\ell.
\end{equation*}
It follows that
\begin{equation}\label{eq:covariant-Hcurl-bound}
	\norm{(D\br)^{\mathsf T}\widehat\bw}_{\bH(\curl,\widehatD)}
	\leq C\norm{D\br}_{L^\infty(\widehatD;\IC^{3\times3})}
	\norm{\widehat\bw}_{\bH^1(\widehatD)},
\end{equation}
with a constant $C>0$ independent of $\br \in W^{1,\infty}(\widehatD;\IC^3)$.
Next, we have that
\begin{equation*}
	\begin{aligned}
		\snorm{\dual{\gamma_{\top,\widehatGamma}\!\left(
				(D\br)^{\mathsf T}\mathsf Z_{\widehatD}\widehat{\bphi}\right)}
			{\widehat\bj}_{\widehatY\times\widehatX}}
		\leq
		C
		\norm{\widehat{\bphi}}_{\bH^{\half}(\widehatGamma)}
		\norm{\widehat\bj}_{\widehatX},
	\end{aligned}
\end{equation*}
where $C>0$ is independent of $\br \in \mathfrak{T}_\delta$.
The map
\begin{equation*}
	W^{1,\infty}(\widehatD;\IC^3)
	\to\mathscr L_{\normalfont\text{bl}}
	\bigl(\widehatX\times\bH^{\half}(\widehatGamma),\IC\bigr):
	\br\mapsto\left((\widehat\bj,\widehat{\bphi})\mapsto
	\dual{\gamma_{\top,\widehatGamma}\!\left((D\br)^{\mathsf T}
	\mathsf Z_{\widehatD}\widehat{\bphi}\right)}
	{\widehat\bj}_{\widehatY\times\widehatX}\right)
\end{equation*}
is continuous and complex linear by the preceding estimates, hence
holomorphic, and bounded on each $\frakT_\delta$.
Applying Lemma~\ref{lem:vector-variational-solution} to
\eqref{eq:coefficient-reference-pairing}, we conclude that there exists
$\delta>0$ such that the map in \eqref{eq:coefficient-holomorphy}
admits a bounded holomorphic extension to $\frakT_\delta$.

Next, for each $\br\in\frakT$, we use \eqref{eq:transported-surface-divergence}
and \eqref{eq:transported-scalar-duality}. Applying Lemma~\ref{lem:piola-div} gives, for
$\widehat\phi\in H^{\half}(\widehatGamma)$,
\begin{equation}\label{eq:div-coefficient-reference-pairing}
	\dual{\widehat{\div}_\br\widehat\bj}{\widehat\phi}_
	{H^{-\half}(\widehatGamma)\times H^{\half}(\widehatGamma),\br}
	=\dual{\divGh\widehat\bj}{\widehat\phi}_
	{H^{-\half}(\widehatGamma)\times H^{\half}(\widehatGamma)}.
\end{equation}
The right-hand side of \eqref{eq:div-coefficient-reference-pairing} is
independent of $\br$, hence bounded and holomorphic on $\frakT_\varepsilon$
for each $\varepsilon>0$. Recalling \cite[Lemma~2.8]{HDPS26Laplace}, we conclude that there exists
$\delta>0$ such that the map in \eqref{eq:div-coefficient-holomorphy}
admits a bounded holomorphic extension to $\frakT_\delta$.
\end{proof}

\section{Shape Holomorphy}
\label{sec:main-results}

In this section, we prove the main result of this work
(Theorem~\ref{thm:shape-holomorphy}), which establishes the holomorphic
dependence of the transported EFIE operator on the domain transformation.
In Subsection~\ref{sec:physical-efie}, the variational EFIE is formulated
on Lipschitz surfaces, and in Subsection~\ref{sec:single-layer}, the shape
holomorphy of the Helmholtz single layer operator is recalled. In Subsection~\ref{sec:efie-factorization},
we combine this result with the transport from Section~\ref{sec:transport}
to prove shape holomorphy of the EFIE operator with bounds explicit in
the wavenumber. In Subsection~\ref{sec:efie-solution}, we use this result
to establish holomorphy of the shape-to-solution map, assuming that the EFIE
is well-posed throughout the collection of admissible shapes and that the
transported data depend holomorphically on the domain transformation.

\subsection{The Electric Field Integral Equation on Lipschitz Surfaces}
\label{sec:physical-efie}
For each $\br\in\frakT$, we consider the following scattering problem.
Let $\einc$ be an incident electric field satisfying $\curl\curl\einc-\kappa^2\einc=0$ in
$\IR^3\setminus\overline{\D_\br}$, with tangential trace
$\gamma_{\top,\Gamma_\br}\einc\in Y_\br$.
We seek a scattered field
$\escat\in\bH_{\mathrm{loc}}(\curl,\IR^3\setminus\D_\br)$ satisfying
\begin{equation}\label{eq:maxwell-scattering-problem}
	\begin{aligned}
		\curl\curl\escat-\kappa^2\escat
		&=0
		&&\text{in }\IR^3\setminus\overline{\D_\br},\\
		\gamma_{\top,\Gamma_\br}\escat
		&=-\gamma_{\top,\Gamma_\br}\einc
		&&\text{on }\Gamma_\br,
	\end{aligned}
\end{equation}
together with the Silver--M\"uller radiation condition
\begin{equation}\label{eq:silver-muller-radiation}
	\norm{(\curl\escat)(\x)\times\frac{\x}{\norm{\x}}
		-\imath\kappa\escat(\x)}
	=o\bigl(\norm{\x}^{-1}\bigr),
	\quad\text{as }\norm{\x}\to\infty,
\end{equation}
uniformly in direction. 

Let $G^{(\kappa)}:\IR^3\setminus\{0\}\to\IC$ be the outgoing fundamental
solution of the Helmholtz operator
\begin{equation*}
	G^{(\kappa)}(\x)\coloneqq
	\frac{e^{\imath\kappa\norm{\x}}}{4\pi\norm{\x}},
	\quad\x\in\IR^3\setminus\{0\}.
\end{equation*}

The single layer potential is defined, for each $\br\in\frakT$ and
smooth densities $\mu$, as
\begin{equation*}
	(\mathsf S_{\Gamma_\br}^{(\kappa)}\mu)(\x)
	\coloneqq
	\int_{\Gamma_\br}
	G^{(\kappa)}(\x-\y)\mu(\y)\ds_\y,
	\quad
	\x\in\IR^3\setminus\Gamma_\br.
\end{equation*}
As is well established (cf.\ \cite[Chapters~6 and~9]{mclean2000strongly}), one has
\begin{equation*}
	\mathsf S_{\Gamma_\br}^{(\kappa)}\in
	\mathscr{L}\bigl(H^{-\half}(\Gamma_\br),H^1(\D_\br)\bigr).
\end{equation*}
We define the single layer boundary integral operator by
\begin{equation*}
	\V_{\Gamma_\br}^{(\kappa)}
	\coloneqq
	\gamma_{\Gamma_\br}\mathsf S_{\Gamma_\br}^{(\kappa)}
	\in
	\mathscr{L}\bigl(
	H^{-\half}(\Gamma_\br),
	H^{\half}(\Gamma_\br)
	\bigr).
\end{equation*}
We use the same notation for the componentwise action of
$\mathsf S_{\Gamma_\br}^{(\kappa)}$ and $\V_{\Gamma_\br}^{(\kappa)}$
on vector-valued densities.

Assume that $\einc$ extends into $\D_\br$ as a Maxwell solution.
With the Stratton--Chu representation formula
\cite[Section~4.2, equation~(37)]{BuffaCostabelSchwab2002}, we may represent
the scattered field in $\IR^3\setminus\overline{\D_\br}$ as
$\escat=\cE^{(\kappa)}_{\Gamma_\br}\bj$, where the electric potential is given by
\begin{equation}\label{eq:electric-potential}
	\cE^{(\kappa)}_{\Gamma_\br}\bj(\x)
	\coloneqq
	\imath\kappa\,
	\mathsf S^{(\kappa)}_{\Gamma_\br}
	(\iota_{\top,\Gamma_\br}\bj)(\x)
	-
	\frac{1}{\imath\kappa}\nabla_{\x}
	\mathsf S^{(\kappa)}_{\Gamma_\br}
	(\div_{\Gamma_\br}\bj)(\x),
	\quad\x\in\IR^3\setminus\Gamma_\br.
\end{equation}
Here $\bj\in X_\br$ is the unknown boundary density,
which corresponds to the tangential trace of the total magnetic field.
The electric potential satisfies the Maxwell equation in
\eqref{eq:maxwell-scattering-problem} and the radiation condition
\eqref{eq:silver-muller-radiation} for every $\bj\in X_\br$.
The interior and exterior tangential traces of the electric potential
agree and belong to $Y_\br$
\cite[Section~3.3]{BHPS2003Maxwell}.
Substituting this representation into the boundary condition in
\eqref{eq:maxwell-scattering-problem}, we obtain the variational EFIE
\cite{BuffaCostabelSchwab2002,BHPS2003Maxwell}: for each $\br\in\frakT$,
find $\bj_\br\in X_\br$ such that
\begin{equation*}\label{eq:EFIE_Problem}
	\mathsf a^{(\kappa)}_{\Gamma_\br}(\bj_\br,\bv)
	=
	-\dual{\gamma_{\top,\Gamma_\br}\einc}{\bv}_{Y_\br\times X_\br},
	\quad\forall\bv\in X_\br,
\end{equation*}
where the bilinear form
$\mathsf a^{(\kappa)}_{\Gamma_\br}:X_\br\times X_\br\to\IC$ is defined by
\begin{equation*}
	\mathsf a^{(\kappa)}_{\Gamma_\br}(\bj,\bv)
	\coloneqq\dual{\gamma_{\top,\Gamma_\br}\cE^{(\kappa)}_{\Gamma_\br}\bj}{\bv}_{Y_\br\times X_\br},
	\quad \bj,\bv\in X_\br.
\end{equation*}
According to \cite[Section~4.2]{BuffaCostabelSchwab2002}, one has,
for all $\bj,\bv\in X_\br$,
\begin{equation}\label{eq:efie-form}
	\begin{aligned}
		\mathsf a^{(\kappa)}_{\Gamma_\br}(\bj,\bv)
		&=
		\imath\kappa
		\dual{
			\iota_{\top,\Gamma_\br}\bv
		}{
			\V^{(\kappa)}_{\Gamma_\br}
			\iota_{\top,\Gamma_\br}\bj
		}_{\bH^{-\half}(\Gamma_\br)
			\times\bH^{\half}(\Gamma_\br)} \\
		&\quad{}+
		\frac{1}{\imath\kappa}
		\dual{
			\div_{\Gamma_\br}\bv
		}{
			\V^{(\kappa)}_{\Gamma_\br}
			\div_{\Gamma_\br}\bj
		}_{H^{-\half}(\Gamma_\br)\times H^{\half}(\Gamma_\br)}.
	\end{aligned}
\end{equation}

The mapping properties of the single layer operator, together with
those of the intrinsic injection and the surface divergence operator
(Lemma~\ref{lem:weak-surface-divergence}), imply that the bilinear form
in \eqref{eq:efie-form} is bounded.
For each $\br\in\frakT$, it therefore induces the operator
$\mathsf A^{(\kappa)}_{\Gamma_\br}\in\mathscr{L}(X_\br,Y_\br)$ by
\begin{equation*}
	\dual{\mathsf A^{(\kappa)}_{\Gamma_\br}\bj}{\bv}_{Y_\br\times X_\br}
	\coloneqq\mathsf a^{(\kappa)}_{\Gamma_\br}(\bj,\bv),
	\quad
	\forall \bj,\bv\in X_\br.
\end{equation*}

\subsection{Shape Holomorphy of the Single Layer Operator}
\label{sec:single-layer}

The results from Section~\ref{sec:transport} require only
Assumption~\ref{ass:geometry}.  To apply the holomorphy results for the scalar single layer operator
from \cite{HDPS26Helmholtz,HDPS26Laplace}, we also require extensions of the
transformations to a common hold-all domain.

\begin{assumption}[{\cite[Assumption~3.1]{HDPS26Laplace}}]
\label{ass:hold-all}
There exists a bounded Lipschitz domain $\cO\subset\IR^3$
such that:
\begin{itemize}
    \item[\normalfont{(i)}]
    $\overline{\widehatD},\overline{\D_\br}\subset\cO$ for each
    $\br\in\frakT$.
    \item[\normalfont{(ii)}]
    Each $\br\in\frakT$ admits an extension, also denoted by
    $\br$, which is a bi-Lipschitz bijection from $\cO$ onto itself.
    The collection of these extensions, denoted by
    $\frakT_{\cO}$, is a compact subset of
    $W^{1,\infty}(\cO;\IR^3)$.
    \item[\normalfont{(iii)}]
    The extensions from item (ii) become the identity in a vicinity of
    $\partial\cO$, i.e., there exists $\varepsilon>0$ such that, for
    each $\br\in\frakT_{\cO}$, it holds that
\begin{equation*}
	\br(\widehat\x)=\widehat\x
	\quad
	\text{for each }
	\widehat\x\in\cO\cap
	\bigcup_{\widehat\y\in\partial\cO}B(\widehat\y,\varepsilon).
\end{equation*}
\end{itemize}
\end{assumption}

For each $\br\in\frakT$, we define the transported single layer operator
\begin{equation}\label{eq:transported-single-layer}
	\widehat{\V}^{(\kappa)}_\br
	\coloneqq
	\tau_\br\V^{(\kappa)}_{\Gamma_\br}\tau_\br^{-1}
	\in
	\mathscr{L}\bigl(
	H^{-\half}(\widehatGamma),
	H^{\half}(\widehatGamma)
	\bigr),
\end{equation}
where, for each $\br \in \mathfrak{T}$, the pullback operator $\tau_\br$
has been introduced in \eqref{eq:scalar-pullback}. 
We also write $\widehat{\V}^{(\kappa)}_\br$ for the componentwise
action of this operator on vector-valued densities.

We recall the following result 
concerning frequency-explicit shape holomorphy 
of the single layer BIO as stated in  \cite[Theorem~3.7]{HDPS26Helmholtz},
which, in turn, extends the result for the Laplace operator established in
\cite[Theorem~3.8]{HDPS26Laplace}.
\begin{theorem}
\label{thm:helmholtz-single-layer}
Let $\mathfrak T$ be a set of domain transformations
satisfying Assumptions~\ref{ass:geometry}
and~\ref{ass:hold-all}.
Then there exists $\delta_0>0$ such that,
for every $\kappa>0$ and $0<\delta<\delta_0$, the map
\begin{equation*}
	\mathcal{V}^{(\kappa)}:
	\mathfrak{T}
	\rightarrow
	\mathscr{L}
	\left(
	H^{-\frac{1}{2}}(\widehat\Gamma),
	H^{\frac{1}{2}}(\widehat\Gamma)
	\right):
	\br
	\mapsto
	\widehat{\mathsf V}^{(\kappa)}_\br
\end{equation*}
admits a bounded holomorphic extension to $\mathfrak T_\delta$.
Moreover, there exists $C_2=C_2(\mathfrak T,\delta_0)>0$ such that, for every
$\kappa_0>0$, there exists $C_1=C_1(\mathfrak T,\delta_0,\kappa_0)>0$
satisfying
\begin{equation}\label{eq:single-layer-frequency-bound}
	\sup_{\br\in\mathfrak T_\delta}
	\norm{\widehat{\mathsf V}^{(\kappa)}_\br}_
	{\mathscr L(H^{-\frac{1}{2}}(\widehat\Gamma),H^{\frac{1}{2}}(\widehat\Gamma))}
	\leq C_1(1+\kappa)
	e^{\frac{4\delta}{\pi R}\left(\kappa\delta_0C_2+2\log(1+\kappa)\right)}
\end{equation}
for all $\kappa\geq\kappa_0$, $0<R<\delta_0$, and
$0<\delta<\min\{R,\delta_0-R\}$.
\end{theorem}

\subsection{Shape Holomorphy of the EFIE}
\label{sec:efie-factorization}
In the following, we use the surface contravariant Piola transform introduced in
\cref{eq:surface-piola-lifting} to transport the EFIE introduced in
\cref{eq:EFIE_Problem}, originally set in $\Gamma_\br$, $\br \in \mathfrak{T}$
to the fixed reference boundary $\widehat{\Gamma}$, while preserving the 
tangential properties of the involved fields. 
To this end, for each $\br\in\frakT$, we define the transported bilinear form
$\widehat{\mathsf a}^{(\kappa)}_\br:\widehatX\times\widehatX\to\IC$ as
\begin{equation}\label{eq:transported-form}
	\widehat{\mathsf a}^{(\kappa)}_\br(\widehat\bj,\widehat\bv)
	\coloneqq
	\mathsf a^{(\kappa)}_{\Gamma_\br}
	\left(
	\mathsf P_\br^{-1}\widehat\bj,\mathsf P_\br^{-1}\widehat\bv
	\right),
	\quad \forall\widehat\bj,\widehat\bv\in\widehatX.
\end{equation}
The boundedness of $\mathsf a^{(\kappa)}_{\Gamma_\br}$ and \Cref{thm:piola} imply that $\widehat{\mathsf a}^{(\kappa)}_\br \in \mathscr{L}_{\normalfont\text{bl}}(\widehat{X} \times \widehat{X},\mathbb{C})$. 
Furthermore, for each $\br \in \mathfrak{T}$, 
recalling \cref{eq:efie-form}, we obtain
\begin{equation}
	\begin{aligned}
		\widehat{\mathsf a}^{(\kappa)}_\br(\widehat\bj,\widehat\bv)
		=
		&
		\imath\kappa
		\dual{
			\iota_{\top,\Gamma_\br} \mathsf{P}_\br^{-1}\widehat\bv
		}{
			\V^{(\kappa)}_{\Gamma_\br}
			\iota_{\top,\Gamma_\br} \mathsf{P}_\br^{-1}\widehat\bj
		}_{\bH^{-\half}(\Gamma_\br)
			\times\bH^{\half}(\Gamma_\br)} \\
		&
		+
		\frac{1}{\imath\kappa}
		\dual{
			\div_{\Gamma_\br} \mathsf{P}_\br^{-1}\widehat\bv
		}{
			\V^{(\kappa)}_{\Gamma_\br}
			\div_{\Gamma_\br} \mathsf{P}_\br^{-1}\widehat\bj
		}_{H^{-\half}(\Gamma_\br)\times H^{\half}(\Gamma_\br)},
		\quad
		\forall \widehat\bj,\widehat\bv\in\widehatX.
	\end{aligned}
\end{equation}
Next, for each $\br \in \mathfrak{T}$, recalling
\cref{eq:transported-scalar-duality} and \eqref{eq:transported-scalar-duality-vector}, 
together with \eqref{eq:coefficient-identification} and \eqref{eq:transported-surface-divergence}, we obtain
\begin{equation}\label{eq:fixed-efie-factorization}
	\begin{aligned}
		\widehat{\mathsf a}^{(\kappa)}_\br(\widehat\bj,\widehat\bv)
		&=
		\imath\kappa
		\dual{
			\widehat{\iota}_{\top,\br}\widehat\bv
		}{
			\widehat{\V}^{(\kappa)}_\br
			\widehat{\iota}_{\top,\br}\widehat\bj
		}_{\bH^{-\half}(\widehatGamma)\times
			\bH^{\half}(\widehatGamma),\br} \\
		&\quad{}+
		\frac{1}{\imath\kappa}
		\dual{
			\widehat{\div}_\br\widehat\bv
		}{
			\widehat{\V}^{(\kappa)}_\br
			\widehat{\div}_\br\widehat\bj
		}_{H^{-\half}(\widehatGamma)\times
			H^{\half}(\widehatGamma),\br},
		\quad \forall \widehat\bj,\widehat\bv\in\widehatX,
	\end{aligned}
\end{equation}
where $\widehat{\iota}_{\top,\br}$ and $\widehat{\div}_\br$ are defined in
\eqref{eq:coefficient-identification} and
\eqref{eq:transported-surface-divergence}, respectively.

For each $\br \in \mathfrak{T}$, we define the transported operator
\begin{equation}\label{eq:transported_EFIE}
	\widehat{\mathsf A}^{(\kappa)}_\br
	\coloneqq \mathsf Q_\br\mathsf A^{(\kappa)}_{\Gamma_\br}\mathsf P_\br^{-1}
	\in
	\mathscr L(\widehatX,\widehatY).
\end{equation}
The application of Lemma~\ref{lem:piola-duality} yields
\begin{equation}\label{eq:EFIE_op_bl_EFIE}
	\dual{
		\widehat{\mathsf A}^{(\kappa)}_\br \widehat\bj
	}{\widehat\bv}_{\widehatY\times\widehatX}
	=
	\widehat{\mathsf a}^{(\kappa)}_\br(\widehat\bj,\widehat\bv),
	\quad
	\forall
	\widehat\bj,\widehat\bv\in\widehatX.
\end{equation}

The next result gives a bounded holomorphic extension of the
operator $\widehat{\mathsf A}^{(\kappa)}_\br$ introduced in
\cref{eq:transported_EFIE} to a complex neighborhood of the form
$\mathfrak{T}_\delta$, for some $\delta>0$.

\begin{theorem}
\label{thm:shape-holomorphy}
Let $\frakT$ be a set of domain transformations satisfying
Assumptions~\ref{ass:geometry} and~\ref{ass:hold-all}.
Then there exists $\delta_0=\delta_0(\frakT)>0$ such that, for every
$\kappa>0$ and $0<\delta<\delta_0$, the map
\begin{equation}\label{eq:map_EFIE_A_hat}
	\cA^{(\kappa)}:
	\frakT\rightarrow\mathscr{L}(\widehatX,\widehatY):
	\br\mapsto\widehat{\mathsf A}^{(\kappa)}_\br
\end{equation}
admits a bounded holomorphic extension to $\frakT_\delta$.
Moreover, there exists $C_2=C_2(\frakT,\delta_0)>0$ such that, for every
$\kappa_0>0$, there exists $C_1=C_1(\frakT,\delta_0,\kappa_0)>0$
satisfying
\begin{equation}\label{eq:efie-frequency-bound}
	\sup_{\br\in\frakT_\delta}
	\norm{\widehat{\mathsf A}^{(\kappa)}_\br}_{\mathscr{L}(\widehatX,\widehatY)}
	\leq C_1(1+\kappa)^2
	e^{\frac{4\delta}{\pi R}\left(\kappa\delta_0C_2+2\log(1+\kappa)\right)}
\end{equation}
for all $\kappa\geq\kappa_0$, $0<R<\delta_0$, and
$0<\delta<\min\{R,\delta_0-R\}$.
\end{theorem}

\begin{proof}
It follows from item (ii) in \Cref{lem:transported-duality,cor:transported-vector-duality}, together with \cref{prop:coefficient-identification,thm:helmholtz-single-layer}, that there exists $\delta_0 = \delta_0(\mathfrak{T})>0$, depending only on $\mathfrak{T}$ and independent of $\kappa$, 
such that for each $\delta\in (0,\delta_0)$ the map
\begin{equation*}
	\mathfrak{T} \rightarrow \mathscr{L}_{\text{bl}}(\widehatX \times \widehatX;\mathbb{C}):
	\br \mapsto \widehat{\mathsf a}^{(\kappa)}_\br
\end{equation*}
admits a bounded holomorphic extension to $\mathfrak{T}_\delta$.
Consequently, in view of \eqref{eq:EFIE_op_bl_EFIE}, we may conclude that there
exists a possibly different $\delta_0$, still independent of $\kappa$, such that
for each $\delta \in (0,\delta_0)$ the map \eqref{eq:map_EFIE_A_hat}
admits a bounded holomorphic extension to $\mathfrak{T}_\delta$.

Next, it follows from \eqref{eq:fixed-efie-factorization} that, for each
$\br\in\frakT_{\delta}$,
\begin{equation*}
	\begin{aligned}
		\norm{\widehat{\mathsf A}^{(\kappa)}_\br}_{\mathscr{L}(\widehatX,\widehatY)}
		\leq C_1
		&\left(
		\kappa
		\norm{\widehat{\iota}_{\top,\br}}^{2}_{\mathscr{L}(\widehatX,\bH^{-\half}(\widehatGamma))}
		\right. \\
		&\left.
		{}+\frac{1}{\kappa}
		\norm{\widehat{\div}_\br}^{2}_{\mathscr{L}(\widehatX,H^{-\half}(\widehatGamma))}
		\right)
		\norm{\widehat{\V}^{(\kappa)}_\br}_{\mathscr{L}(H^{-\half}(\widehatGamma),H^{\half}(\widehatGamma))},
	\end{aligned}
\end{equation*}
where $C_1$ is as in \cref{eq:transported-duality-bound}.
Recalling item (i) in
\Cref{lem:transported-duality,cor:transported-vector-duality},
\cref{prop:coefficient-identification,thm:helmholtz-single-layer}, and
\cref{eq:single-layer-frequency-bound}, we obtain
\cref{eq:efie-frequency-bound} with a different constant $C_1>0$ which is
still independent of $\kappa\geq\kappa_0$ but depending on $\kappa_0$.
\end{proof}

\subsection{Shape Holomorphy of the EFIE Solution}
\label{sec:efie-solution}

We now use Theorem~\ref{thm:shape-holomorphy} to pass from shape holomorphy
of the EFIE operator to that of its solution. For a fixed wavenumber, we
assume that the EFIE is well-posed on each admissible boundary. The following
result establishes a bounded holomorphic extension of the inverse operator
and, consequently, of the solution map whenever the transported right-hand
side has the same regularity. 

\begin{corollary}[Shape holomorphy of the EFIE solution on closed surfaces]
\label{cor:efie-solution}
Let $\frakT$ satisfy
Assumptions~\ref{ass:geometry} and~\ref{ass:hold-all}.
Suppose that, for each $\br\in\frakT$, $\kappa^2$ is not an interior
electric Maxwell eigenvalue of $\D_\br$.
That is, the only $\bE\in\bH_0(\curl,\D_\br)$ satisfying
\begin{equation*}
	\begin{aligned}
		\int_{\D_\br}\curl\bE(\x)\cdot \overline{\curl\bV(\x)}\,\dd\x
		=
		\kappa^2\int_{\D_\br}\bE(\x)\cdot \overline{\bV(\x)}\,\dd\x,
		\quad\forall\bV\in\bH_0(\curl,\D_\br),
	\end{aligned}
\end{equation*}
is $\bE=0$.
Then the following assertions hold.
\begin{itemize}
 \item[\normalfont{(i)}] There exists $\delta>0$ such that the map
\begin{equation}
	\frakT\to\mathscr L(\widehatY,\widehatX):
	\br\mapsto
	(\widehat{\mathsf A}^{(\kappa)}_\br)^{-1}
\end{equation}
 admits a bounded holomorphic extension to $\frakT_\delta$.
 \item[\normalfont{(ii)}] For each $\br\in\frakT$, let
 $\bff_\br\in Y_\br$ and assume
 that $\frakT\ni\br\mapsto\widehat{\bff}_\br\coloneqq\mathsf Q_\br\bff_\br\in\widehatY$
 admits a bounded holomorphic extension to $\frakT_{\delta_f}$
 for some $\delta_f>0$. For each $\br\in\frakT$, let
 $\bj_\br\in X_\br$ be the unique solution to
\begin{equation}
	\mathsf A^{(\kappa)}_{\Gamma_\br}\bj_\br
	=\bff_\br.
\end{equation}
 Then there exists $\delta\in(0,\delta_f]$
 such that the map
\begin{equation}
	\frakT\to\widehatX:\br\mapsto
	\widehat\bj_\br\coloneqq\mathsf P_\br\bj_\br
	=(\widehat{\mathsf A}^{(\kappa)}_\br)^{-1}
	\widehat{\bff}_\br
\end{equation}
 admits a bounded holomorphic extension to $\frakT_\delta$.
\end{itemize}
\end{corollary}

\begin{proof}
Under the stated eigenvalue condition,
$\mathsf A^{(\kappa)}_{\Gamma_\br}:X_\br\to Y_\br$ is a bounded
isomorphism for every $\br\in\frakT$
\cite[Corollary~5.5]{BHPS2003Maxwell}.
By Proposition~\ref{thm:piola} and the definition of the transported EFIE
operator in \eqref{eq:transported_EFIE},
$\widehat{\mathsf A}^{(\kappa)}_\br$ is a bounded isomorphism for every
$\br\in\frakT$. Theorem~\ref{thm:shape-holomorphy}, compactness of
$\frakT$, and the perturbation estimate in
\cite[Proposition~4.20, item~(i)]{HS21} yield $\delta>0$ such that the
holomorphic operator extension is invertible throughout $\frakT_\delta$,
with uniformly bounded inverses. Holomorphy of inversion
\cite[Proposition~4.20, item~(ii)]{HS21} proves item~(i).

After decreasing $\delta$ so that $\delta\leq\delta_f$, item~(ii) follows
by applying the bounded bilinear evaluation map to the inverse and the
transported right-hand side
\cite[Proposition~4.20, item~(iii)]{HS21}.
\end{proof}

For a plane wave $\einc(\x)
=\bp\exp(\imath\kappa\boldsymbol d\cdot \x)$ with
$\boldsymbol d\in\mathbb S^2$ and $\bp\in\IC^3$ satisfying
$\bp\cdot \boldsymbol d=0$, set
$\bff_\br=-\gamma_{\top,\Gamma_\br}\einc$.
The transported data satisfy
\begin{equation*}
	\dual{\widehat{\bff}_\br}{\widehat\bv}_{\widehatY\times\widehatX}
	=-\dual{\widehat{\iota}_{\top,\br}\widehat\bv}{\btau_\br\einc}_
	{\bH^{-\half}(\widehatGamma)\times\bH^{\half}(\widehatGamma),\br},
	\quad \widehat\bv\in\widehatX.
\end{equation*}
This identity follows from \eqref{eq:coefficient-identification},
the transported duality \eqref{eq:transported-scalar-duality},
the characterization \eqref{eq:iota-top-characterization}, and
Lemma~\ref{lem:piola-duality} with $\bseta=\bff_\br$ and
$\bj=\mathsf P_\br^{-1}\widehat\bv$.
The exponential map $\br\mapsto\exp(\imath\kappa\boldsymbol d\cdot\br)$ is holomorphic.
The boundedness of the Dirichlet trace,
Corollary~\ref{cor:transported-vector-duality}, and
Proposition~\ref{prop:coefficient-identification} therefore give
a bounded holomorphic extension of $\widehat{\bff}_\br$ to a possibly
smaller $\frakT_\delta$.

\section{Parametric Holomorphy}
\label{sec:efie-parametric}

We now apply the shape holomorphy results of Section~\ref{sec:main-results}
to affine domain transformations depending on finitely or countably many
parameters. Under suitable summability assumptions on the affine terms,
we establish $(\boldsymbol{b},p,\varepsilon)$-holomorphy of the
parameter-to-operator map and, under the hypotheses of
Corollary~\ref{cor:efie-solution}, of the parameter-to-solution map.
Throughout this section, we set
$\mathbb{U}\coloneqq[-1,1]^{\mathbb N}\subset\ell^\infty(\mathbb N;\IR)$
and equip $\mathbb{U}$ with the product topology. By Tychonoff's theorem,
$\mathbb{U}$ is compact in this topology.

\subsection{Affine-Parametric Domain Transformations}

Let $\cO\subset\IR^3$ be the bounded Lipschitz hold-all domain from
Assumption~\ref{ass:hold-all}. Given
$\br_0\in W^{1,\infty}(\cO;\IR^3)$ and
$(\boldsymbol\psi_j)_{j\geq 1}\subset W^{1,\infty}(\cO;\IR^3)$,
we define, for each $\boldsymbol y=(y_j)_{j\geq 1}\in\mathbb{U}$,
the affine-parametric domain transformation
\begin{equation}\label{eq:efie-affine-shapes}
	\br_{\boldsymbol y}
	\coloneqq
	\br_0+\sum_{j\geq 1}y_j\boldsymbol\psi_j
	\quad\text{and}\quad
	\boldsymbol{b}\coloneqq(b_j)_{j\geq 1},
	\quad
	b_j\coloneqq
	\norm{\boldsymbol\psi_j}_{W^{1,\infty}(\cO;\IR^3)}.
\end{equation}
We assume that $\boldsymbol{b}\in\ell^p(\mathbb{N})$ for some
$p\in(0,1)$ and that the set
\begin{equation}\label{eq:efie-parametric-family}
	\mathfrak{T}
	\coloneqq
	\left\{\br_{\boldsymbol y}|_{\widehatD}:\boldsymbol y\in \mathbb{U}\right\}
\end{equation}
satisfies Assumptions~\ref{ass:geometry} and
\ref{ass:hold-all}. When used as an element of $\mathfrak{T}$,
$\br_{\boldsymbol y}$ denotes its restriction to $\widehatD$.
Since $\boldsymbol{b}\in\ell^p(\mathbb{N})\subset\ell^1(\mathbb{N})$,
the series in \eqref{eq:efie-affine-shapes} converges absolutely in
$W^{1,\infty}(\cO;\IR^3)$, uniformly with respect to
$\boldsymbol y\in\mathbb{U}$. The resulting map
$\boldsymbol y\mapsto\br_{\boldsymbol y}$ is continuous from
$\mathbb{U}$, equipped with the product topology, into
$W^{1,\infty}(\cO;\IR^3)$. Compactness of $\mathbb{U}$ therefore implies
compactness of the family on $\cO$ and of its restrictions to
$\widehatD$ \cite[Section~5.2]{CSZ18}.

\paragraph{Affine-parametric perturbations of a Lipschitz domain}
To construct a family satisfying these assumptions, let $\widehatD$
be a bounded Lipschitz domain with connected boundary and
$\overline{\widehatD}\subset\cO$.
Multiplying a bounded linear extension by a fixed smooth cutoff equal
to one near $\overline{\widehatD}$, we obtain an extension operator
$\mathscr E:W^{1,\infty}(\widehatD;\IR^3)\to W^{1,\infty}(\IR^3;\IR^3)$
with images supported in $K\Subset\cO$.
For $\boldsymbol\varphi_j\in W^{1,\infty}(\widehatD;\IR^3)$, set
$a_j\coloneqq\norm{\boldsymbol\varphi_j}_{W^{1,\infty}(\widehatD;\IR^3)}$
and assume $(a_j)_{j\geq1}\in\ell^p(\mathbb N)$.
For $\theta>0$ and $\boldsymbol y\in \mathbb{U}$, define
\begin{equation*}
	\br_{\boldsymbol y}(\widehat\x)
	\coloneqq\widehat\x+
	\theta\sum_{j\geq1}y_j(\mathscr E\boldsymbol\varphi_j)(\widehat\x),
	\qquad\widehat\x\in\IR^3.
\end{equation*}
The resulting family has the form \eqref{eq:efie-affine-shapes} with
$\br_0=\operatorname{Id}$,
$\boldsymbol\psi_j=\theta(\mathscr E\boldsymbol\varphi_j)|_{\cO}$, and
$b_j\leq \theta C_{\mathscr E}a_j$, where $C_{\mathscr E}$ denotes the
operator norm of $\mathscr E$. We verify
Assumptions~\ref{ass:geometry} and~\ref{ass:hold-all}
for sufficiently small $\theta>0$. Compactness follows from the preceding
argument. According to \cite[item (i) in Lemma~B.1]{HDPS26Laplace},
there exists $\delta_{\mathrm L}>0$ such that every
$\br\in W^{1,\infty}(\widehatD;\IR^3)$ satisfying
$\norm{\br-\operatorname{Id}}_{W^{1,\infty}(\widehatD;\IR^3)}<\delta_{\mathrm L}$
is a bi-Lipschitz bijection from $\widehatD$ onto the bounded Lipschitz
domain $\D_\br=\br(\widehatD)$. We choose $\theta>0$ so that
\begin{equation*}
	\sup_{\boldsymbol y\in \mathbb{U}}
	\norm{\br_{\boldsymbol y}-\operatorname{Id}}
	_{W^{1,\infty}(\IR^3;\IR^3)}
	\leq \theta C_{\mathscr E}\sum_{j\geq1}a_j
	\eqqcolon q<\min\{\delta_{\mathrm L},1/2\}.
\end{equation*}
Since $\operatorname{Lip}(\br_{\boldsymbol y}-\operatorname{Id})\leq q<1$,
the contraction mapping theorem gives global bi-Lipschitz bijections
$\br_{\boldsymbol y}$, so
$\Gamma_{\br_{\boldsymbol y}}=\br_{\boldsymbol y}(\widehatGamma)$ is connected.
Moreover,
$\norm{D(\br_{\boldsymbol y}-\operatorname{Id})}_{L^\infty(\IR^3;\IR^{3\times3})}\leq q<1$,
which gives
\begin{equation*}
	\inf_{\boldsymbol y\in \mathbb{U}}
	\operatorname*{ess\,inf}_{\widehat\x\in\widehatD}
	J_{\br_{\boldsymbol y}}(\widehat\x)
	\geq(1-q)^3>0.
\end{equation*}
Since each map equals the identity outside $K$, the transformations
are the identity in a common neighborhood of $\partial\cO$ of width
$d_K\coloneqq\operatorname{dist}(K,\partial\cO)>0$.
Global bijectivity then gives $\br_{\boldsymbol y}(\cO)=\cO$ and
$\overline{\D_{\br_{\boldsymbol y}}}\subset\cO$.

\subsection{Parametric Holomorphy of the EFIE Operator and Solution}

We recall the notion of parametric holomorphy introduced in
\cite[Definition~2.1]{CCS15}, following the formulation in
\cite[Definition~2.6]{dolz2024parametric}.
For $\rho>1$, let
\begin{equation}
	\mathcal{E}_\rho
	\coloneqq
	\left\{\frac{\zeta+\zeta^{-1}}{2}:
	1\leq\snorm{\zeta}\leq\rho\right\}\subset\mathbb{C}
\end{equation}
be the Bernstein ellipse with foci at $-1$ and $1$. For a sequence
$\boldsymbol{\rho}=(\rho_j)_{j\geq 1}$ with $\rho_j>1$, we define the
polyellipse $\mathcal{E}_{\boldsymbol{\rho}}\coloneqq
\bigotimes_{j\geq 1}\mathcal{E}_{\rho_j}$.

\begin{definition}\label{def:bpe_holomorphy}
Let $X$ be a complex Banach space, let $p\in(0,1)$ and
$\varepsilon>0$, and let
$\boldsymbol{b}=(b_j)_{j\geq 1}\in\ell^p(\mathbb{N})$ be a
nonnegative sequence. A map
$\Phi:\mathbb{U}\rightarrow X$ is called
$(\boldsymbol{b},p,\varepsilon)$-holomorphic if the following assertions
hold.
\begin{itemize}
    \item[\normalfont{(i)}]
    The map $\Phi$ is uniformly bounded on $\mathbb{U}$.
    \item[\normalfont{(ii)}]
    For every sequence $\boldsymbol{\rho}=(\rho_j)_{j\geq 1}$ satisfying
\begin{equation}\label{eq:b_epsilon_admissibility}
	\rho_j>1
	\quad\text{and}\quad
	\sum_{j\geq 1}(\rho_j-1)b_j\leq\varepsilon,
\end{equation}
    the map $\Phi$ admits an extension, also denoted by $\Phi$, which is
    holomorphic with respect to each variable $z_j$ on
    $\mathcal{O}_{\boldsymbol{\rho}}
    \coloneqq\bigotimes_{j\geq 1}\mathcal{O}_{\rho_j}$, where
    $\mathcal{O}_{\rho_j}\coloneqq
    \left\{z\in\mathbb{C}:\operatorname{dist}(z,[-1,1])<\rho_j-1\right\}$.
    \item[\normalfont{(iii)}]
    There exists a constant $C_\varepsilon>0$, independent of
    $\boldsymbol{\rho}$, such that
\begin{equation}\label{eq:bpe_uniform_bound}
	\sup_{\boldsymbol z\in\mathcal{E}_{\boldsymbol{\rho}}}
	\norm{\Phi(\boldsymbol z)}_X\leq C_\varepsilon.
\end{equation}
\end{itemize}
\end{definition}

We now use the affine representation \eqref{eq:efie-affine-shapes}
to pass from shape holomorphy to
$(\boldsymbol{b},p,\varepsilon)$-holomorphy.

\begin{proposition}\label{prop:shape_to_parametric_holomorphy}
Let $X$ be a complex Banach space and let
$\mathcal{A}:\mathfrak{T}\rightarrow X$ admit a bounded holomorphic
extension to $\mathfrak{T}_\delta$ for some $\delta>0$. Let
$\mathfrak{T}$ be given by \eqref{eq:efie-parametric-family} with
$\boldsymbol{b}\in\ell^p(\mathbb{N})$ for some $p\in(0,1)$. Then, for
every $\varepsilon\in(0,\delta)$, the map
\begin{equation}
	\Phi:\mathbb{U}\rightarrow X:
	\boldsymbol y\mapsto\mathcal{A}(\br_{\boldsymbol y})
\end{equation}
is $(\boldsymbol{b},p,\varepsilon)$-holomorphic.
\end{proposition}

\begin{proof}
Let $\boldsymbol{\rho}$ satisfy
\eqref{eq:b_epsilon_admissibility} and let
$\boldsymbol z\in\mathcal{O}_{\boldsymbol{\rho}}$. Set
$\br_{\boldsymbol z}\coloneqq\br_0+\sum_{j\geq1}z_j\boldsymbol\psi_j$.
The series converges in $W^{1,\infty}(\cO;\IC^3)$ by
\eqref{eq:b_epsilon_admissibility} and
$\boldsymbol{b}\in\ell^p(\mathbb{N})\subset\ell^1(\mathbb{N})$.
For each $j\geq1$, choose
$y_j\in[-1,1]$ such that
$\snorm{z_j-y_j}<\rho_j-1$. Then $\boldsymbol y=(y_j)_{j\geq 1}\in \mathbb{U}$ and
\begin{equation}
	\norm{\br_{\boldsymbol z}-\br_{\boldsymbol y}}_{W^{1,\infty}(\cO;\IC^3)}
	\leq
	\sum_{j\geq 1}\snorm{z_j-y_j}b_j
	\leq\sum_{j\geq1}(\rho_j-1)b_j
	\leq\varepsilon<\delta.
\end{equation}
Consequently, $\br_{\boldsymbol z}|_{\widehatD}\in\mathfrak{T}_\delta$.
Since the affine map
$\boldsymbol z\mapsto\br_{\boldsymbol z}|_{\widehatD}$ is holomorphic
with respect to each variable, composition with the bounded
holomorphic extension of $\mathcal{A}$ gives item~(ii) of
Definition~\ref{def:bpe_holomorphy}, as in
\cite[Section~5.2]{CSZ18}.
The boundedness of this extension gives item~(i) and, since
$\mathcal{E}_{\rho_j}\subset\mathcal{O}_{\rho_j}$ for each $j\geq1$,
item~(iii).
\end{proof}

We now apply Proposition~\ref{prop:shape_to_parametric_holomorphy}
to the transported EFIE operator and its solution.

\begin{corollary}\label{cor:efie-parametric}
Assume that the family in \eqref{eq:efie-parametric-family} satisfies
Assumptions~\ref{ass:geometry} and~\ref{ass:hold-all}, with
$\boldsymbol b\in\ell^p(\mathbb N)$ for some $p\in(0,1)$.
For each $\kappa>0$, there exists $\varepsilon>0$ such that the
parameter-to-operator map
\begin{equation*}
	\mathbb{U}\to\mathscr L(\widehatX,\widehatY):
	\boldsymbol y\mapsto\widehat{\mathsf A}^{(\kappa)}_{\br_{\boldsymbol y}}
\end{equation*}
is $(\boldsymbol b,p,\varepsilon)$-holomorphic.
If the hypotheses of Corollary~\ref{cor:efie-solution} hold, then the
parameter-to-solution map
\begin{equation*}
	\mathbb{U}\to\widehatX:\boldsymbol y\mapsto
	\widehat\bj_{\br_{\boldsymbol y}}
\end{equation*}
is $(\boldsymbol b,p,\varepsilon)$-holomorphic, possibly with a smaller
$\varepsilon$. 
\end{corollary}

\begin{proof}
The assertions follow by applying
Proposition~\ref{prop:shape_to_parametric_holomorphy} to the bounded
holomorphic extensions provided by
Theorem~\ref{thm:shape-holomorphy} and
Corollary~\ref{cor:efie-solution}, respectively.
\end{proof}

\section{Concluding Remarks}
\label{sec:concluding-remarks}

We have established shape holomorphy of the EFIE operator
on the boundaries of bounded Lipschitz domains in three dimensions,
including polyhedra. This result extends the analyses in
\cite{Henriquez2026C11,EscapilJerez2026EFIE} beyond the $\mathscr C^{1,1}$ setting.
At the core of the proof is the construction of surface Piola
transformations through volume liftings and the curl-conforming pullback.
We combine these transformations with the transported duality of
\cite{HDPS26Laplace} and the single layer holomorphy result of
\cite{HDPS26Helmholtz} to obtain the extension on the Maxwell trace spaces
of the fixed reference boundary.

In addition, we provide bounds for the transported EFIE operator that are
explicit in the wavenumber and the size of the domain of holomorphic
extension. This domain can be chosen independently of $\kappa>0$.
The estimate in \eqref{eq:efie-frequency-bound} describes how the operator
norm grows with these two quantities. Choosing $\delta$ proportional to
$\kappa^{-1}$ keeps the exponential factor uniformly bounded and yields
an operator norm bound of order $(1+\kappa)^2$ for $\kappa\geq\kappa_0>0$.

These results also allow us to establish holomorphic dependence of the
EFIE solution on the geometry. At a fixed wavenumber, assuming that the
EFIE is well posed for every admissible domain transformation,
Corollary~\ref{cor:efie-solution} gives a bounded holomorphic extension
of the inverse to a possibly smaller complex neighborhood. Applying this
inverse to holomorphic transported data then yields shape holomorphy of
the transported boundary density. In particular, the required regularity
of the data holds for plane-wave incidence.

Under the summability assumptions of Section~\ref{sec:efie-parametric},
composing the shape-holomorphic operator and solution maps with affine
domain transformations gives $(\boldsymbol b,p,\varepsilon)$-holomorphy
of the corresponding parameter maps. The approximation results of
\cite{CCS15} then yield sparse polynomial approximations with convergence
rates independent of the number of parameters. This regularity also
provides a basis for reduced order models that approximate the EFIE
boundary density over families of geometries
\cite{WKHH25,DH2026,PintoHenriquez2025ReducedBasis}.

Current work encompasses the extension of this framework to the EFIE on
Lipschitz screens, extending the corresponding analysis in the
$\mathscr C^{1,1}$ setting \cite{Henriquez2026C11}, and to the full Maxwell
Calder\'on projector, with a view to treating transmission problems.

\section*{Declaration on the use of generative AI}
ChatGPT (OpenAI) assisted with the weak normal trace argument in Lemma 2.1, parts of the proof of Lemma 2.9, the recovery of the classical surface Piola formulas in Remark 2.10, literature searches, and manuscript revision. The author verified the mathematical arguments and takes responsibility for all content.

\raggedbottom

\end{document}